\documentclass{commat}

\usepackage{longtable}
\usepackage{stmaryrd}
\usepackage{cleveref}
\usepackage{tikz}    
\usetikzlibrary{arrows}

\newtheorem*{conj*}{Conjecture}

\theoremstyle{definition}
\newtheorem{exm}{Example}

\newcommand{\FF}{\mathbb{F}}

\newcommand{\NN}{\mathbb{N}}

\newcommand{\A}{\mathsf{A}}

\newcommand{\af}{\alpha}

\newcommand{\Dt}{\Delta}
\newcommand{\e}{\epsilon}
\newcommand{\s}{\sigma}

\newcommand{\ad}{\mathrm{ad}}

\newcommand{\tl}{\tilde}

\newcommand{\sst}{\subseteq}
 
\newcommand{\hder}[1]{\operatorname{\mathrm{HDer}}{#1}}
\newcommand{\ihder}[1]{\operatorname{\mathrm{IHDer}}{#1}}

\newcommand{\hd}[1]{\{{#1}_n\}_{n=0}^\infty}

\newcommand{\cS}{\mathcal{S}}
\newcommand{\cT}{\mathcal{T}}

\newcommand{\impl}{\Rightarrow}

\newcommand{\lb}{\lambda}

 \usepackage[symbol]{footmisc}

\title[Non-associative algebraic structures:
classification and structure]{%
   Non-associative algebraic structures:
classification and structure\footnote{	The work was supported by FCT UIDB/MAT/00212/2020, UIDP/MAT/00212/2020, and 2022.02474.PTDC.}
    }

\author{%
    Ivan Kaygorodov
    }

\authorinfo[%
I. Kaygorodov
    ]{
    CMA-UBI, Universidade da Beira Interior, Covilh\~{a}, Portugal}{%
    kaygorodov.ivan@gmail.com
    }

\abstract{%
    These are detailed notes for a~lecture on “Non-associative Algebraic Structures:
Classification and Structure” which I
 presented as part of my Agregação em Matemática e Applicações
(University of Beira Interior, Covilhã, Portugal, 13-14/03/2023).
    }



\begin{document}

\VOLUME{32}
\YEAR{2024}
\NUMBER{3}
\firstpage{1}
\DOI{https://doi.org/10.46298/cm.11419}

 \tableofcontents
 
\section{Introduction}
In this lecture, we aim to summarize our research interests and achievements as well as motivate what drives our work.
The present work is based on various papers
written together with
Abror Khudoyberdiyev,
Alexandre Pozhidaev,
Aloberdi Sattarov,
Amir Fern\'andez Ouaridi,
Antonio Jesús Calderón,
Artem Lopatin,
Bakhrom Omirov,
Bruno Ferreira,
Crislaine Kuster,
Doston Jumaniyozov,
Elisabete Barreiro,
Farukh Mashurov,
Feng Wei,
Henrique Guzzo,
Gulkhayo Solijanova,
Iqboljon Karimjanov,
Ilya Gorshkov,
Inomjon Yuldashev,
Isabel Hern\'{a}ndez,
Isamiddin Rakhimov,
Ivan Shestakov,
Jobir Adashev,
José María Sánchez,
Karimbergan Kudaybergenov,
Luisa Camacho,
Manat Mustafa,
Manuel Ladra,
María Alejandra Alvarez,
Mikhail Ignatyev,
Mohamed Salim,
Mykola Khrypchenko,
Nurlan Ismailov,
Pasha Zusmanovich,
Patrícia Beites,
Paulo Saraiva,
Pilar P\'{a}ez-Guill\'{a}n,
Renato Fehlberg Júnior,
Samuel Lopes,
Sharifah Kartini Said Husain,
Thiago Castilho de Mello,
Ualbai Umirbaev,
Vasily Voronin,
Viktor Lopatkin,
Yury Popov and
Yury Volkov
~\cite{IKP21,KPV21,CFK191,KKY,ACK21,fkk22,IKM21,KK20cdgeo,KK20cdalg,KM20,CKKK20,KKS20,KIS20,JKK21geo,AK21,AKKS21,AKK21,JKK21,com5,anti6,nilp4,kv17,kv16,kpv,kppv,klp,kpv19,kkp19geo,kkk18,jkk19,ikv18,ikm19,ikv17,gkp,gkk19,gkk19alg,fkkv, cfk19,maria,ack,geo21,CFK21,kv19,klp20, Kay1,kantorII,kpp19,kayvo,klp18mjm,KZ21,cklos20,ccko21,KLP18inv,KKW19,FK21alt,FGK21alt,KP,kay14mz, kay14st,bkp19,kaypopov16, KP16,FKK,fkl21,KK21p,CKKK22,KKL21,KLP20ce,ksu18,ckls20,K17,bks19,bcks19,cks19,KP19,KPS19}.

The main direction of our research is concentrated on non-associative algebraic structures, their classifications, and the study of their properties.
Namely,
we work with many types of non-associative algebras and superalgebras,
$n$-ary algebras and also algebraic structures with two binary multiplications
(known as algebras of Poisson type).
We consider small-dimensional algebras, simple and semisimple algebras, PI-algebras, and free algebras.

\subsection{Non-associative algebras}
Unless otherwise noted, all vector spaces and algebras are considered over the complex field.
In our work, \textquotedblleft algebra\textquotedblright \,
will mean \textquotedblleft non-associative algebra\textquotedblright, i.e.,
a vector space $\A$ equipped with a~bilinear product which need not be associative or unital.
In the case of associative algebras, we work with algebras defined by the following
identity: $ (xy)z=x(yz).$
On the other hand, in non-associative cases, there are algebras defined by some special conditions (not identities!),
such as axial algebras, evolution algebras, composition algebras, conservative algebras, etc., and also algebras defined by polynomial identities.
It seems that the first example of non-associative algebras defined by some family of polynomial identities are Lie algebras.
All these algebras are closely related to associative algebras.
For example, each associative algebra under the commutator product
$[x,y]:=xy-yx$ gives a~Lie algebra
(the variety of algebras satisfying this property is called Lie-admissible algebras).

Commutative associative algebras and Lie algebras admit many generalizations, which have triggered some interest.
So, the varieties of
Jordan algebras,
alternative algebras,
right alternative algebras,
assosymmetric algebras,
bicommutative algebras,
Novikov algebras,
left-symmetric algebras,
right commutative algebras,
$\mathfrak{CD}$-algebras,
weakly associative algebras,
noncommutative Jordan algebras or
terminal algebras,
are generalizations of commutative associative algebras.
On the other hand,
the varieties of
Malcev,
symmetric Leibniz,
Leibniz,
binary Lie,
$\mathfrak{CD}$-algebras, and also
noncommutative Jordan algebras and
terminal algebras
can be seen as generalizations of Lie algebras.
Some types of algebras appeared in another way:
Zinbiel algebras are Koszul dual to Leibniz algebras,
Tortkara algebras are defined by identities that are found in Zinbiel algebras under the commutator multiplication;
mock Lie algebras are commutative analogs of Lie algebras;
and dual mock Lie algebras are Koszul dual to mock Lie algebras.
There are many other relations between cited varieties, such as the following.
Lie algebras and Jordan algebras are related by
Kantor–Koecher–Tits construction.
Alternative and noncommutative Jordan algebras under the symmetric product
$ x \bullet y=xy+yx$ give Jordan algebras,
assosymetric algebras under the symmetric product give commutative $\mathfrak{CD}$-algebras;
Novikov, bicommutative, assosymmetric, and left-symmetric algebras are Lie-admissible algebras (i.e., under the commutator product they give Lie algebras) and so on.

Below, we provide the list of identities that defined the varieties of algebras that are under consideration in this manuscript.

\begin{longtable}
{|lcl|}
\hline
$(-1,1)$- &:&
$(xy)y = xy^2,\ (x,y,z) + (z,x,y) + (y,z,x) = 0$\\\hline

Alternative &:&
$ (x,y,z)=-(y,x,z), \ (x,y,z)=-(x,z,y)$\\\hline
Antiassociative &:&
$ (xy)z=-x(yz)$\\\hline
Assosymmetric &:&
$ (x,y,z)=(y,x,z), \ (x,y,z)=(x,z,y)$ \\\hline
Bicommutative &:&
$ x(yz)=y(xz), \ (xy)z=(xz)y$ \\\hline
Binary Lie&:& $xy=-yx, ((xy)y)x=((xy)x)y$\\\hline
$\mathfrak{CD}$- &:& $((xy)a)b-((xy)b)a=((xa)b-(xb)a)y+x((ya)b-(yb)a)$,\\
& & $ (a(xy))b-a((xy)b)=((ax)b-a(xb))y+x((ay)b-a(yb))$,\\
& & $ a(b(xy))-b(a(xy))=(a(bx)-b(ax))y+x(a(by)-b(ay))$\\\hline
Dual mock-Lie &:& $xy=-yx, (xy)z=-x(yz)$\\\hline
Jordan & :&
$ xy=yx, \ (x^2 y)x=x^2 (yx)$\\ \hline
Left-symmetric &:&
$ (x,y,z)=(y,x,z)$ \\\hline
Leibniz&:& $(xy)z=(xz)y+x(yz)$\\ \hline
Malcev &:& $xy=-yx, \ (xy)(xz)+(y(xz))x+$ \\ && \multicolumn{1}{c|}{$+((xz)x)y=((xy)z+(yz)x+(zx)y)x$} \\\hline
Noncommutative Jordan &:&
$(xy)x=x(yx), \ (x^2 y)x=x^2 (yx)$\\\hline
Novikov &:&
$ (xy)z=(xz)y, \ (x,y,z)=(y,x,z)$ \\\hline
Right alternative &:&
$ (x,y,z)=-(x,z,y)$ \\\hline
Right commutative &:&
$ (xy)z=(xz)y $\\\hline
Symmetric Leibniz&:& $x(yz)=(xy)z+y(xz), \ (xy)z=(xz)y+x(yz)$\\\hline

Terminal &:&
$3(b(a(xy)-(ax)y-x(ay)) - a((bx)y) + (a(bx))y+$\\
&& \multicolumn{1}{c|}{ $(bx)(ay) -a(x(by))+(ax)(by)+x(a(by)))=$}\\
&& \multicolumn{1}{r|}{$= - (2ab+ba)(xy)+((2ab+ba)x)y + x((2ab+ba)y)$}\\\hline
Tortkara &:& $xy=-yx, \ (xy)(zy)=((xy)z+(yz)x+(zx)y)y $\\\hline

Weakly associative &:&
$(xy)z-x(yz)+(yz)x-y(zx) = (yx)z- y(xz)$\\\hline

Zinbiel &:& $(xy)z=x(yz+zy)$\\ \hline
\end{longtable}

\textbf{Superalgebras.}
Let $\mathbb{F}$ be a~field and let $G$ be the Grassmann algebra over $\mathbb{F}$ given by the generators
$1 $, $\xi_1$, $\ldots $, $\xi_n $, $\ldots$ and the defining relations $\xi_i^2 =0$ and $\xi_i\xi_j=-\xi_j\xi_i $.
The elements $1$, $\xi_{i_1} \xi_{i_2} \ldots \xi_{i_k} $, $i_1< i_2 < \ldots < i_k$, form a~basis of the algebra $G$ over $\mathbb{F}$. Denote by $G_0$ and $G_1$ the subspaces spanned by the products
of even and odd lengths, respectively; then $G$ can be represented as the direct sum of these subspaces,
$G = G_0 \oplus G_1 $.
Here the relations
\[
G_iG_j \subseteq G_{i+j (mod \ 2)}, \ i,j = 0, 1,
\]
hold.
In other words, $G$ is a~$\mathbb{Z}_2$-graded algebra (or a~superalgebra) over $\mathbb{F}$.
Suppose now that $\A = \A_0 \oplus \A_1$ is an arbitrary superalgebra over $\mathbb{F}$. Consider the tensor product $G \otimes \A$ of $\mathbb{F}$-algebras.
The subalgebra
\[
G(\A) = G_0 \otimes \A_0 + G_1 \otimes \A_1
\]
of $G \otimes \A$ is referred to as the Grassmann envelope of the superalgebra $\A$.
Let $\Omega$ be a~variety of algebras over $\mathbb{F}$.
A superalgebra $\A = \A_0 \oplus \A_1$ is referred to as an
$\Omega$-superalgebra if its Grassmann envelope $G(\A)$ is an algebra in $\Omega$.
In particular, $\A = \A_0 \oplus \A_1$ is
referred to as a~(noncommutative) Jordan superalgebra if its Grassmann envelope $G(\A)$ is a~(noncommutative) Jordan algebra.

\subsection{Poisson algebras}
	The notion of a~Poisson bracket has its origin in the works of S.D. Poisson on celestial mechanics at the beginning of the XIX century. Since then, Poisson algebras
	have appeared in several areas of mathematics, such as:
	Poisson manifolds~\cite{L1},
	algebraic geometry~\cite{BLLM,GK04},
	noncommutative geometry~\cite{v08},
	operads~\cite{MR},
	quantization theory~\cite{Hue90,Kon03},
	quantum
	groups~\cite{Dr87},
	and classical and quantum mechanics.
	The study of Poisson algebras also led to other algebraic structures, such as
	generic Poisson algebras~\cite{ksu18,KSML},
	algebras of Jordan brackets and generalized Poisson algebras~\cite{KacKant07,ZheKa11,marsheze01,KingMcc92},
	Gerstenhaber algebras~\cite{ksy},
	Novikov-Poisson algebras~\cite{Zakharov},
	Malcev-Poisson-Jordan algebras,
	transposed Poisson algebras~\cite{bai20,fkl21},
	$n$-ary Poisson algebras~\cite{ck16}, etc.

\begin{definition}
	A Poisson algebra is a~vector space $\A$ with two bilinear operations $\cdot$ and $\{ \cdot,\cdot\}$, such that $(\A,\cdot)$ is a~commutative associative algebra, $(\A,\{ \cdot, \cdot \})$ is a~Lie algebra and
\[
	\{x,y\cdot z\}=\{x,y\}\cdot z+y\cdot\{x,z\}.
\]

\end{definition}

Poisson algebras are at the intersection of Lie and commutative associative algebras.
They inherit many common properties of these algebras,
and they can be used to solve some well-known problems in the theory of Lie or associative algebras.
So,
	Poisson algebras have been used to prove the Nagata conjecture about wild automorphisms of the polynomial ring with three generators~\cite{su04}
	and to describe simple Jordan~\cite{MZ01}
	and noncommutative Jordan~\cite{PS13,PS19} algebras and superalgebras.
    Namely, from each simple Poisson algebra, it is possible to obtain simple Jordan superalgebra by a~process named as \textquotedblleft Kantor double\textquotedblright; but on the other hand, each noncommutative Jordan superalgebra under the commutator product gives a~structure of Poisson bracket.
	The systematic study of free Poisson algebras began in the paper by Shestakov~\cite{s93}, where he constructed a~basis for the free Poisson algebra.
	Later, in a~series papers by Makar-Limanov and Umirbaev (see~\cite{mlu09,mltu09} and references therein),
	many analogs of classical results for free Poisson algebras were obtained.
	For example, they proved that there are no wild automorphisms in the free Poisson algebra with two generators, proved the freedom theorem (The Freiheitssatz), and described universal multiplicative envelopings of the free Poisson fields.
	The systematic study of noncommutative Poisson algebra structures began in the paper by Kubo~\cite{kubo96}.
	He obtained a~description of all the Poisson structures on the full and upper triangular matrix algebras, which was later generalized to prime associative algebras in~\cite{fale98}.
	Namely, it was proved in~\cite{fale98} that any Poisson bracket on a~prime noncommutative associative algebra is the commutator bracket multiplied by an element of the extended centroid of the algebra.
	On the other hand, in his next paper, Kubo studied
	noncommutative Poisson algebra structures on affine Kac-Moody algebras.
	The investigation of Poisson structures on associative algebras continued in some papers of
	Yao, Ye and Zhang~\cite{yyz07}; Mroczyńska, Jaworska-Pastuszak, and Pogorzały~\cite{jpp20,mp18},
	where Poisson structures on finite-dimensional path algebras
	and on canonical algebras were studied.
	Crawley-Boevey introduced a~noncommutative Poisson structure, called an $H_0$-Poisson structure, on the $0$-th cyclic homology of a noncommutative associative algebra~\cite{cb11}
	and showed that an $H_0$-Poisson structure can be induced on the affine moduli space of (semisimple) representations of an associative algebra from a~suitable Lie algebra structure on the $0$-th Hochschild homology of the algebra~\cite{cb11}.
	The derived noncommutative Poisson bracket on Koszul Calabi-Yau algebras has been studied in~\cite{ceey17}.
	Van den Bergh introduced double Poisson algebras in~\cite{v08}, and
	Van de Weyer described all the double Poisson structures on finite-dimensional semisimple algebras.
	Recently, the notion of a~noncommutative Poisson bialgebra appeared in~\cite{lbs20}.

\subsection{$n$-Ary algebras}
Following Kurosh,
an $\Omega$-algebra over a~field $\mathbb F$ is a~linear space over $\mathbb F$ equipped with a~system of multilinear algebraic operations $\Omega =\{ \omega_i: |\omega_i|=n_i, n_i \in \mathbb N\} $,
where $|\omega_i|$ denotes the arity of $\omega_i$.
The present notion of $\Omega$-algebras is very wide.
It includes associative and non-associative algebras (with one binary multiplication);
dialgebras and Poisson algebras and their generalizations (with two binary multiplications);
trialgebras (with three binary multiplications);
quadri-algebras (with four binary multiplications),
ennea-algebras (with nine binary multiplications);
pluriassociative algebras (with $2\gamma$ binary multiplications);
Jordan triple disystems and comtrans algebras (with two ternary multiplications);
Akivis algebras (with one binary and one ternary multiplications);
Bol algebras,
Lie-Yamaguti algebras,
Poisson $n$-Lie algebras,
Sabinin algebras and so on.
A particular case of $\Omega$-algebras, that will be interesting for us is
$n$-ary algebras.

Let $\A$ be an $\Omega$-algebra equipped with one $n$-ary operation
$[\cdot, \ldots, \cdot]: \A \times \ldots \times \A \to \A$,
which is $n$-linear.
This operation is $n$-ary
commutative if,
\[
[x_1, \ldots, x_n]= [x_{\sigma(1)}, \ldots, x_{\sigma(n)}], \mbox{ for all } \sigma \in \mathbb S_n.
\]

There are notions of partial commutativity and partial associativity for $n$-ary algebras.
A generalization of a variety of binary algebras defined by a~family of polynomial identities to the $n$-ary case is not an easy task.
So, there are many types of generalizations of associative and Lie algebras to $n$-ary case.
For example, there are partially associative and totally associative $n$-ary algebras, which give associative algebras in the binary case.
In the case of Lie algebras,
the known generalizations are the following:
reduced Lie ternary algebras introduced by Pojidaev and
$n$-Lie algebras introduced by Filippov in~\cite{Fil}.
 Recently, the theory of Filippov algebras has attracted much attention due to its close connection with Nambu mechanics, proposed by Nambu as a~generalization of classical Hamiltonian mechanics.
 In his article, Filippov introduced the notion of a~$n$-Lie algebra and proved some of its properties. He also obtained the classification of anticommutative $n$-ary algebras of dimension $n$ and $n+1$.

Let us give the definition of $n$-Lie (Filippov) algebras.
Let $\A$ be an algebra equipped with one $n$-ary operation
$[\cdot, \ldots, \cdot]: \A \times \ldots \times \A \to \A$,
which is $n$-linear.
This algebra is $n$-Lie if,
\[
[x_1, \ldots, x_n]= (-1)^{\sigma} [x_{\sigma(1)}, \ldots, x_{\sigma(n)}], \mbox{ for all } \sigma \in \mathbb S_n;
\]
 \[
  [[x_1, \ldots, x_n], y_2, \ldots, y_n]= \sum\limits_{i=1}^n [x_1, \ldots, x_{i-1}, [x_i, y_2, \ldots, y_n], x_{i+1}, \ldots, x_n].
 \]
 The first identity gives an $n$-ary analog of anticommutative identity, and the second gives the $n$-ary version of the Jacobi identity.
 An algebra $\A$ is called $n$-ary Leibniz algebra if it satisfies only the second identity from the identities mentioned above.

\section{Classifications of non-associative algebras}

The present part is based on the papers
written together with
Abror Khudoyberdiyev,
Alexandre Pozhidaev,
Aloberdi Sattarov,
Amir Fern\'andez Ouaridi,
Antonio Jesús Calderón,
Crislaine Kuster,
Doston Jumaniyozov,
Farukh Mashurov,
Iqboljon Karimjanov,
Ilya Gorshkov,
Isabel Hern\'{a}ndez,
Isamiddin Rakhimov,
Jobir Adashev,
Luisa Camacho,
Manat Mustafa,
Manuel Ladra,
María Alejandra Alvarez,
Mikhail Ignatyev,
Mohamed Salim,
Mykola Khrypchenko,
Nurlan Ismailov,
Pilar P\'{a}ez-Guill\'{a}n,
Samuel Lopes,
Sharifah Kartini Said Husain,
Thiago Castilho de Mello,
Vasily Voronin,
Yury Popov and
Yury Volkov~\cite{IKP21,KPV21,CFK191,ACK21,fkk22,IKM21,KK20cdgeo,KK20cdalg,KM20,CKKK20,KKS20,KIS20,JKK21geo,AK21,AKKS21,AKK21,JKK21,com5,anti6,nilp4,kv17,kv16,kpv,kppv,klp20,kpv19,kkp19geo,kkk18,jkk19,ikv18,ikm19,ikv17,CKKK22,KKL21,KLP20ce,gkp,gkk19,gkk19alg,fkkv, cfk19,maria,ack,geo21,CFK21,kv19,ckls20}.
\subsection{The algebraic classification of algebras}
The algebraic classification (up to isomorphism) of algebras of dimension $n$ from a certain variety
defined by a certain family of polynomial identities is a~classic problem in the theory of non-associative algebras.
There are many results related to the algebraic classification of small-dimensional algebras in many varieties of
non-associative algebras~\cite{ ck13, degr3, usefi1, degr2, ha16, kv16}.
So, algebraic classifications of
$2$-dimensional algebras~\cite{kv16,petersson},
$3$-dimensional evolution algebras~\cite{ccsmv},
$3$-dimensional anticommutative algebras~\cite{japan},
$4$-dimensional division algebras~\cite{Ernst,erik},
 have been given.

Our method for classifying nilpotent commutative algebras is based on the calculation of central extensions of nilpotent algebras of smaller dimensions from the same variety (first, this method has been developed by Skjelbred and Sund for the Lie algebra case in~\cite{ss78}).

Our big program in the classification of complex small-dimensional nilpotent algebras results in three big classifications of the following algebras:
\begin{itemize}
    \item $4$-dimensional nilpotent algebras~\cite{nilp4};
    \item $5$-dimensional commutative nilpotent algebras~\cite{com5};
    \item $6$-dimensional anticommutative nilpotent algebras~\cite{anti6}.
\end{itemize}

On the other hand, during the realization of our classification program,
we classified small-dimensional algebras in the following varieties:
\begin{itemize}
    \item[I.] $4$-dimensional nilpotent:
    \begin{itemize}
        \item[---] assosymmetric algebras~\cite{ikm19};
        \item[---] bicommutative algebras~\cite{kpv19};
        \item[---] $\mathfrak{CD}$-algebras~\cite{KK20cdalg};
        \item[---] commutative algebras~\cite{fkkv};
        \item[---] left-symmetric algebras~\cite{AKKS21};
        \item[---] noncommutative Jordan algebras~\cite{jkk19};
        \item[---] Novikov algebras~\cite{kkk18};
        \item[---] right commutative algebras~\cite{AKK21};
        \item[---] right alternative algebras~\cite{IKM21};
        \item[---] terminal algebras~\cite{kkp19geo};
        \item[---] weakly associative algebras~\cite{AK21};

          \end{itemize}
    \item[II.] $5$-dimensional nilpotent:
    \begin{itemize}
        \item[---] antiassociative algebras~\cite{fkk22};
        \item[---] anticommutative algebras~\cite{fkkv};
        \item[---] commutative associative algebras~\cite{KIS20};
        \item[---] commutative $\mathfrak{CD}$-algebras~\cite{JKK21};
        \item[---] symmetric Leibniz algebras~\cite{AK21};
        \item[---] one-generated assosymmetric algebras~\cite{KM20};
        \item[---] one-generated bicommutative algebras~\cite{KPV21};
        \item[---] one-generated Novikov algebras~\cite{CKKK20};
        \item[---] one-generated terminal algebras~\cite{KKS20};

  \end{itemize}
    \item[III.] $6$-dimensional nilpotent:
    \begin{itemize}
        \item[---] anticommutative $\mathfrak{CD}$-algebras~\cite{ack};
        \item[---] binary Lie algebras~\cite{ack};
        \item[---] Tortkara algebras~\cite{gkk19alg};
        \item[---] one-generated assosymmetric algebras~\cite{KM20};
        \item[---] one-generated bicommutative algebras~\cite{KPV21};
        \item[---] one-generated Novikov algebras~\cite{CKKK20}.

         \end{itemize}

\end{itemize}

The classification of non-nilpotent algebras is given by other methods, which were developed in our papers.
Let us now summarize the rest of the results obtained in the algebraic classification of algebras.

\begin{itemize}
    \item[$\bullet$] $2$-dimensional terminal, conservative and rigid algebras~\cite{cfk19};
    \item[$\bullet$] $2$-dimensional algebras with left quasi-units~\cite{kayvo};
    \item[$\bullet$] $2$-dimensional algebras~\cite{kv16};
    \item[$\bullet$] $3$-dimensional anticommutative algebras\cite{ikv18};
    \item[$\bullet$] $8$-dimensional dual Mock Lie algebras~\cite{ckls20};
    \item[$\bullet$] central extensions of $n$-dimensional filiform associative algebras~\cite{KKL21};
    \item[$\bullet$] central extensions of $n$-dimensional filiform Zinbiel algebras~\cite{CKKK22};
    \item[$\bullet$] non-associative central extensions of $n$-dimensional null-filiform associative algebras~\cite{KLP20ce};

    \item[$\bullet$] $n$-dimensional algebras with $(n-2)$-dimensional annihilator~\cite{CFK21};

    \item[$\bullet$] $n$-dimensional anticommutative algebras with $(n-3)$-dimensional annihilator~\cite{CFK191};
    \item[$\bullet$] $n$-dimensional Zinbiel algebras with $(n-3)$-dimensional annihilator~\cite{ACK21}.
\end{itemize}

\subsection{The geometric classification of algebras}

Given a~complex $n$-dimensional vector space $\mathbb V$, the set \[ {\rm Hom}(\mathbb V \otimes \mathbb V,\mathbb V) \cong \mathbb V^* \otimes \mathbb V^* \otimes \mathbb V \]
is a~complex vector space of dimension $n^3$. This space has the structure of the affine variety $\mathbb{C}^{n^3}$. Indeed, if we fix a~basis $\{e_1,\dots,e_n\}$ of $\mathbb V$, then any $\mu\in {\rm Hom}(\mathbb V \otimes \mathbb V,\mathbb V)$ is determined by $n^3$ structure constants $c_{ij}^k \in\mathbb{C}$ such that
$\mu(e_i\otimes e_j)=\sum\limits_{k=1}^n c_{ij}^k e_k$. A subset of ${\rm Hom}(\mathbb V \otimes \mathbb V,\mathbb V)$ is \textit{Zariski-closed} if it is the set of solutions of a~system of polynomial equations in the variables $c_{ij}^k$ ($1\le i,j,k\le n$).

Let $T$ be a set of polynomial identities.
Every algebra structure on $\mathbb V$ satisfying polynomial identities from $T$ forms a~Zariski-closed subset of the variety ${\rm Hom}(\mathbb V \otimes \mathbb V,\mathbb V)$. We denote this subset by $\mathbb{L}(T)$.
The general linear group ${\rm GL}(\mathbb V)$ acts on $\mathbb{L}(T)$ by conjugation:
\[
(g * \mu)(x\otimes y) = g\mu(g^{-1}x\otimes g^{-1}y)
\]
for $x,y\in \mathbb V$, $\mu\in \mathbb{L}(T)\subset {\rm Hom}(\mathbb V \otimes\mathbb V, \mathbb V)$ and $g\in {\rm GL}(\mathbb V)$.
Thus, $\mathbb{L}(T)$ decomposes into ${\rm GL}(\mathbb V)$-orbits that correspond to the isomorphism classes of the algebras.
We shall denote by $O(\mu)$ the orbit of $\mu\in\mathbb{L}(T)$ under the action of ${\rm GL}(\mathbb V)$ and by $\overline{O(\mu)}$ its Zariski closure.

One of the main problems of the geometric classification of a variety of algebras is a description of its irreducible components. In~\cite{gabriel}, Gabriel described the irreducible components
of the variety of $4$-dimensional unital associative algebras, and the variety of $5$-dimensional
unital associative algebras were classified algebraically and geometrically by Mazzola~\cite{maz79,maz80}.
Later, Cibils~\cite{Cibils} considered rigid associative algebras with 2-step nilpotent radical. Goze and
Ancochea-Bermúdez proved that the varieties of 7 and 8-dimensional nilpotent Lie algebras
are reducible~\cite{esp}. The irreducible components of $2$-step nilpotent commutative associative
algebras were described in~\cite{shaf}.
Often, the irreducible components of the variety are determined by the rigid algebras,
although this is not always the case. Indeed, Flanigan showed in~\cite{flan} that the variety
of $3$-dimensional nilpotent associative algebras has an irreducible component which does not
contain any rigid algebras — it is instead defined by the closure of a~union of a~one-parameter
family of algebras.

Let $\A_1$ and $\A_2$ be two $n$-dimensional algebras satisfying the identities from $T$, and let $\mu,\lambda \in \mathbb{L}(T)$ represent $\A_1$ and $\A_2$, respectively.
We say that $\A_1$ \textit{degenerates to} $\A_2$,
and write $\A_1\to \A_2$,
if $\lambda\in\overline{O(\mu)}$.
Note that this implies $\overline{O(\lambda)}\subset\overline{O(\mu)}$. Hence, the definition of degeneration does not depend on the choice of $\mu$ and $\lambda$.
If $\A_1\not\cong \A_2$, then the assertion $\A_1\to \A_2$ is called a~\textit{proper degeneration}.
Following Gorbatsevich~\cite{gorb91}, we say that $\A$ has \textit{level} $m$ if there exists a~chain of proper degenerations of length $m$ starting in $\A$ and there is no such chain of length $m+1$. Also, in~\cite{gorb93} it was introduced the notion of \emph{infinite level} of an algebra $\A$ as the limit of the usual levels of $\A \oplus \mathbb{C}^m$.

Let $\A$ be represented by $\mu\in\mathbb{L}(T)$. Then $\A$ is \textit{rigid} in $\mathbb{L}(T)$ if $O(\mu)$ is an open subset of $\mathbb{L}(T)$.
 Recall that a~subset of a~variety is called \textit{irreducible} if it cannot be represented as a~union of two non-trivial closed subsets.
 A maximal irreducible closed subset of a~variety is called an \textit{irreducible component}.
It is well known that any affine variety can be uniquely represented as a~finite union of its irreducible components. Note that the algebra $\A$ is rigid in $\mathbb{L}(T)$ if and only if $\overline{O(\mu)}$ is an irreducible component of $\mathbb{L}(T)$.

In our work, we discuss the following problems:

\begin{problem}[Geometric classification]
Let $\mathfrak{V}^n$ be a~variety of $n$-dimen\-sio\-nal algebras defined by a~family of identities $T$.
Which are the irreducible components of $\mathfrak{V}^n$?
\end{problem}

\begin{problem}[Degeneration classification]
Let $\mathfrak{V}^n$ be a~variety of $n$-dimen\-sio\-nal algebras defined by a~family of identities $T$.
Is there a~degeneration from $\A_1$ to $\A_2$ for each pair of algebras from $\mathfrak{V}^n$?
\end{problem}

\begin{problem}[Level classification]
Let $\mathfrak{V}^n$ be a~variety of $n$-dimensional algebras defined by a~family of identities $T$.
Which algebras from $\mathfrak{V}^n$ have level $m$?
\end{problem}

 In the following list, we summarize our main results related to the geometric classification of algebras.
 Namely, we gave the geometric classification of the following varieties of algebras:

\begin{itemize}
    \item[I.] $4$-dimensional nilpotent:
     \begin{itemize}
        \item[---] assosymmetric algebras~\cite{ikm19};
        \item[---] bicommutative algebras~\cite{kpv19};
        \item[---] $\mathfrak{CD}$-algebras~\cite{KK20cdgeo};
        \item[---] commutative algebras~\cite{fkkv};
        \item[---] left-symmetric algebras~\cite{AKKS21};
        \item[---] Leibniz algebras~\cite{kppv};
        \item[---] noncommutative Jordan algebras~\cite{jkk19};
        \item[---] Novikov algebras~\cite{kkk18};
        \item[---] right commutative algebras~\cite{AKK21};
        \item[---] right alternative algebras~\cite{IKM21};
        \item[---] terminal algebras~\cite{kkp19geo};
        \item[---] weakly associative algebras~\cite{AK21};

          \end{itemize}
    \item[II.] $5$-dimensional nilpotent:
    \begin{itemize}
        \item[---] antiassociative algebras~\cite{fkk22};
        \item[---] associative algebras~\cite{IKP21};

        \item[---] commutative $\mathfrak{CD}$-algebras~\cite{JKK21geo};
        \item[---] symmetric Leibniz algebras~\cite{AK21};

  \end{itemize}
    \item[III.] $6$-dimensional nilpotent:
    \begin{itemize}
        \item[---] anticommutative $\mathfrak{CD}$-algebras~\cite{ack};
        \item[---] binary Lie algebras~\cite{ack};
        \item[---] Malcev algebras~\cite{kpv};
        \item[---] Tortkara algebras~\cite{gkk19};
     \end{itemize}

    \item[IV.] $4$-dimensional:
    \begin{itemize}
    \item[---] binary Lie algebras~\cite{kpv};
    \item[---] Leibniz algebras~\cite{ikv17};
    \item[---] Zinbiel algebras~\cite{kppv};
     \end{itemize}

    \item[V.] $7$-dimensional:
    \begin{itemize}
         \item[---] dual Mock Lie algebras~\cite{ckls20}.
     \end{itemize}

\end{itemize}

In~\cite{geo21,IKP21} we completely solve the geometric classification problem for
nilpotent and $2$-step nilpotent, commutative
nilpotent and anticommutative nilpotent algebras of arbitrary dimension.

\begin{theorem}[Theorem A,~\cite{geo21}]
For any $n\ge 2$, the variety of all $n$-dimensional nilpotent algebras is irreducible and has dimension $\frac{n(n-1)(n+1)}{3}$.
\end{theorem}

Moreover, we show that the family $\mathcal{R}_n$ given in~\cite[Definition 10]{geo21} is generic in the variety of $n$-dimensional nilpotent algebras and inductively gives an algorithmic procedure to obtain any $n$-dimensional nilpotent algebra as a~degeneration from $\mathcal{R}_n$.

\begin{theorem}[Theorem B,~\cite{geo21}]
For any $n\ge 2$, the variety of all $n$-dimensional commutative nilpotent algebras is irreducible and has dimension $\frac{n(n-1)(n+4)}{6}$.
\end{theorem}

As above, we show that the family $\mathcal{S}_n$ given in~\cite[Definition 15]{geo21} is generic in the variety of $n$-dimensional commutative nilpotent algebras and inductively give an algorithmic procedure to obtain any $n$-dimensional nilpotent commutative algebra as a~degeneration from $\cS_n$.

\begin{theorem}[Theorem C,~\cite{geo21}]
For any $n\ge 2$, the variety of all $n$-dimensional anticommutative nilpotent algebras is irreducible and has dimension $\frac{(n-2)(n^2 +2n+3)}{6}$.
\end{theorem}

We show also that the family $\mathcal{T}_n$ given in~\cite[Definition 33]{geo21}, is generic in the variety of $n$-dimensional anticommutative nilpotent algebras and inductively give an algorithmic procedure to obtain any $n$-dimensional nilpotent anticommutative algebra as a~degeneration from $\cT_n$.

The notion of length for non-associative algebras has been recently introduced in~\cite{guterman}, generalizing the corresponding notion for associative algebras. Using the above result, we show in~\cite[Corollary 39]{geo21} that the length of an arbitrary (i.e., not necessarily nilpotent) $n$-dimensional anticommutative algebra is bounded above by the $n^{\text{th}}$ Fibonacci number, and prove that our bound is sharp.

For $k \le n$ consider the (algebraic) subset $\mathfrak{Nil}_{n,k}^2$ of the variety $\mathfrak{Nil}_n^2$ of 2-step nilpotent $n$-dimensional algebras defined by
\[
\mathfrak{Nil}_{n,k}^2 = \{ \A \in \mathfrak{Nil}_n^2 : \dim \A^2 \le k, \ \dim {\rm ann} \ \A \ge k \}.
\]
It is easy to see that $\mathfrak{Nil}_n^2 = \cup_{k=1}^n \mathfrak{Nil}_{n,k}^2 $. Analogously, for the varieties $\mathfrak{Nil}_n^2 {}^{\mathrm{c}}, \mathfrak{Nil}_n^2 {}^{\mathrm{ac}}$ of commutative and anticommutative 2-step nilpotent algebras we define the subsets $\mathfrak{Nil}_{n,k}^2 {}^{\mathrm{c}}$ and $\mathfrak{Nil}_{n,k}^2 {}^{\mathrm{ac}}$, respectively.

\begin{theorem}[Theorem A,~\cite{IKP21}]
The sets $\mathfrak{Nil}_{n,k}^2$ are irreducible and
\begin{longtable}
{lclll}
$\mathfrak{Nil}_n^2$ &$ =$&$ \bigcup_{k}\mathfrak{Nil}_{n,k}^2 $,& for &$ 1 \le k \le \left\lfloor n + \frac{1 - \sqrt{4n+1}}{2}\right\rfloor$\\
$\mathfrak{Nil}_n^2 {}^{\mathrm{c}}$ &$=$& $\bigcup_{k}\mathfrak{Nil}_{n,k}^2 {}^{\mathrm{c}}$, & for & $1 \le k \le \left\lfloor n + \frac{3 - \sqrt{8n+9}}{2}\right\rfloor$,\\
$\mathfrak{Nil}_n^2 {}^{\mathrm{ac}}$ &$=$ & $\bigcup_{k}\mathfrak{Nil}_{n,k}^2 {}^{\mathrm{ac}}$, & for & $1 + (n + 1) \operatorname{mod} 2 \le k \le \left\lfloor n + \frac{1 - \sqrt{8n+1}}{2}\right\rfloor \text{ for } n \ge 3$.\\
\end{longtable}
Moreover,
\begin{longtable}
{lcl}
$\dim \mathfrak{Nil}_{n,k}^2$ &$=$ &$(n-k)^2 k + (n-k)k$, \\
$\dim \mathfrak{Nil}_{n,k}^2 {}^{\mathrm{c}}$ &$=$ & $\frac{(n-k)(n-k+1)}{2}k + (n-k)k $,\\
$\dim \mathfrak{Nil}_{n,k}^2 {}^{\mathrm{ac}}$ &$=$ & $\frac{(n-k)(n-k-1)}{2}k + (n-k)k $.
\end{longtable}
\end{theorem}


\subsection{Degenerations of algebras}

The present part is dedicated to the results obtained
regarding Problem 2 (Degeneration classification) and
Problem 3 (Level classification).
Typically, the results of degeneration classification are given by construction of the graph of primary degenerations of the algebras from the variety.
So,
Grunewald and O'Halloran constructed the graph of primary degenerations for $5$-dimensional nilpotent Lie algebras~\cite{GRH};
Seeley obtained a~similar graph for $6$-dimensional nilpotent Lie algebras in~\cite{S90}.
The degenerations graphs for $2$-dimensional Jordan algebras and $4$-dimensional nilpotent Jordan algebras are presented in~\cite{esp2,esp3}.
Other graphs of degenerations were constructed by Burde and co-authors
for $4$-dimensional Lie algebras~\cite{BC99},
$2$-dimensional Pre-Lie algebras~\cite{bb09} and
$3$-dimensional Novikov algebras~\cite{bb14}.

In the following list, we summarize our main results related to the degeneration classification of algebras.
 Namely, we constructed graphs of primary degenerations of the following varieties of algebras:

\begin{itemize}
    \item[$\bullet$] $2$-dimensional algebras~\cite{kv16};
    \item[$\bullet$] $3$-dimensional anticommutative algebras~\cite{ikv18};
    \item[$\bullet$] $3$-dimensional Jordan algebras~\cite{gkp};
    \item[$\bullet$] $3$-dimensional Jordan superalgebras~\cite{maria};
    \item[$\bullet$] $3$-dimensional Leibniz algebras~\cite{ikv18};
    \item[$\bullet$] $3$-dimensional nilpotent algebras~\cite{fkkv};
    \item[$\bullet$] $4$-dimensional binary Lie algebras~\cite{kpv};
    \item[$\bullet$] $4$-dimensional nilpotent commutative algebras~\cite{fkkv};
    \item[$\bullet$] $4$-dimensional nilpotent Leibniz algebras~\cite{kppv};
    \item[$\bullet$] $4$-dimensional Zinbiel algebras~\cite{kppv}.
    \item[$\bullet$] $5$-dimensional nilpotent anticommutative algebras~\cite{fkkv};
    \item[$\bullet$] $5$-dimensional nilpotent commutative associative algebras~\cite{klp20}.
    \item[$\bullet$] $6$-dimensional nilpotent Malcev algebras~\cite{kpv};
    \item[$\bullet$] $7$-dimensional dual Mock Lie algebras~\cite{ckls20};
    \item[$\bullet$] $(n+1)$-dimensional $n$-Lie (Filippov) algebras~\cite{kv19}.

\end{itemize}

The first attempt at classification of algebras of level one is given by Gorbatsevich~\cite{gorb91},
but he omitted some main cases.
Later, his classification was completed in a~paper by Khudoyberdiyev and Omirov~\cite{khud13}.
After this,
the classification of algebras of level two for some particular cases,
such that associative algebras, Jordan algebras, Lie algebras, Leibniz algebras and nilpotent algebras has been given in~\cite{khud15,fkrv}.
In~\cite{kv17} we completely solve the level classification problem for
all algebras in the case of level two, generalizing all the cited results of Khudoyberdiyev and co-authors.
Following our ideas on the level classification of algebras, Volkov recently obtained
some results on the classifications of $n$-ary algebras of level one and
on classification of anticommutative algebras of levels $3 $, $4$ and $5$~\cite{wolf1,wolf2}.

\section{Conservative algebras and superalgebras}
The present part is based on the papers
written together with
Alexandre Pozhidaev,
Artem Lopatin,
Renato Fehlberg Júnior,
Yury Popov, and
Yury Volkov~\cite{klp, Kay1,kantorII,kpp19,kayvo}.

\subsection{Kantor product}\label{kanpro}

The idea of obtaining new objects from old ones by using derivative operations has long been
known in algebra~\cite{albert}. In its most general form, the idea was realized by Malcev~\cite{malcev}. Let $M_n$
be an associative algebra of matrices of order $n$ over a~field $\mathbb F $. Assume that some finite collection
$\Lambda=(a_{ij},b_{ij},c_{ij})$ of matrices in $M_n$ is given. Denote by $M^{(\Lambda)}_n$ the algebra defined on a~space of
matrices in $M_n$ with respect to new multiplication
\[
x \cdot _{\Lambda} y =
\sum_{i,j} a_{ij}xb_{ij}yc_{ij}.
\]
It was proved
that every $n$-dimensional algebra over $\mathbb F$ is isomorphic to a~subalgebra of $M^{(\Lambda)}_n$~\cite{malcev}.
Other interesting ways to derive the initial multiplication are
isotopes, homotopes, and mutations~\cite{elduque91, alsaody, pche2}.
The concept of an isotope was introduced by Albert~\cite{albert}.
Let $\A$ and $\A_0$ be algebras with a~common underlying
linear space on which right multiplication operators $R_x$ and $R^{(0)}_x$ are defined (for $\A$ and $\A_0$, resp.).
We say that $\A_0$ and $\A$ are isotopic if there exist invertible linear operators $\phi, \psi, \xi$ such that
$R^{(0)}_x = \phi R_{x\psi} \xi $.
We call $\A_0$ an isotope of $\A$.
Let $\A$ be an arbitrary associative algebra,
and let $p$, $q$ be two fixed elements of $\A $.
Then a
new algebra is derived from $\A$ by using the same vector space structure of $\A$
but defining a~new multiplication
$x * y = x p y - y q x $.
The resulting algebra is called the
$(p, q)$-mutation of the algebra $\A$.

The definition of the Kantor product of multiplications comes from the study of a~certain class of algebras. In 1972, Kantor introduced the class of conservative algebras~\cite{Kan1}, which contains many important classes of algebras (see~\cite{klp}), for example, associative, Lie, Jordan, and Leibniz algebras. To define what will be called the Kantor product, we need to introduce the algebras $U(n)$ (see, for more details,~\cite{Kan2,klp}). Consider the space $U(n)$ of all bilinear multiplications on the $n$-dimensional vector space $V_n$. Now, fix a~vector $u\in V_n$. For $A,B\in U(n)$ (two multiplications) and $x,y\in V_n$, we set
\[
x\ast y=\llbracket A,B \rrbracket (x,y)=A(u,B(x,y))-B(A(u,x),y)-B(x,A(u,y)).
\]

This new multiplication is called the (left) Kantor product of the multiplications $A$ and $B$ (it is also possible to define the right Kantor product). The Kantor product of a~multiplication ``$\cdot$'' by itself will be the Kantor square of ``$\cdot$'':
\[
x\ast y=u\cdot(x\cdot y)-(u\cdot x)\cdot y-x\cdot (u\cdot y).
\]
It is easy to see that the Kantor square of multiplication is a~particular case of the Malcev construction in a~non-associative sense.
On the other hand,
 \begin{itemize}
    \item in the commutative associative case it coincides with a~mutation;
\item in the left commutative and left-symmetric cases it coincides with an isotope.
\end{itemize}
As in~\cite{Kay1}, we will assume that the Kantor product is always the left Kantor product.
In~\cite{Kay1}, we studied the Kantor product and Kantor square of many well-known algebras, for example, associative, (anti)-commutative, Lie, Leibniz, Novikov, Poisson algebras, and dialgebras.

We summarize the main results of the present part in the following theorems.

\begin{theorem}[Theorem 29,~\cite{Kay1}]
Let $(\A, \cdot)$ be a~finite-dimensional associative algebra.
Then $(\A, \cdot)$ is isomorphic to $(\A,\llbracket \cdot,\cdot \rrbracket)$,
if and only if $A$ is a~skew field.

\end{theorem}

\begin{theorem}[Theorem 10,~\cite{Kay1}]
Let $(\A;\cdot)$ be an alternative algebra.
    Then $(\A;\llbracket\cdot, \cdot \rrbracket)$ is a~flexible algebra. Furthermore,
\begin{itemize}
    \item $(\A;\llbracket\cdot, \cdot \rrbracket)$ is an alternative algebra if and only if $\A$ satisfies the identity
    \begin{center}
$(x,u,(x,u,y))=0 $;
\end{center}
    \item $(\A;\llbracket\cdot, \cdot \rrbracket)$ is a~noncommutative Jordan algebra if and only if $\A$
    satisfies the identity
    \begin{center}
$[L_uL_xL_uL_x, R_uR_x] = [L_{xuxu}, R_{ux}] $;
\end{center}
    \item $(\A;\llbracket\cdot, \cdot \rrbracket)$ is a~Jordan algebra if $(\A, \cdot)$ is a~commutative alternative algebra$ $;
    \item $(\textbf{C};\llbracket\cdot, \cdot \rrbracket)$ is alternative for a~Cayley --- Dickson algebra $\textbf{C}$,
    if and only if $u= u_0\cdot 1 $.
    \end{itemize}
    \end{theorem}

\begin{theorem}[Theorem 19,~\cite{Kay1}]
Let $(\A;\cdot,\{\cdot,\cdot\})$ be a~generalized Poisson algebra.
    Then $(\A;\llbracket\{\cdot,\cdot\}, \cdot \rrbracket)$ is an associative-commutative algebra, and
    $(\A;\llbracket\cdot, \{\cdot,\cdot\} \rrbracket)$ is a~Lie algebra.
    \end{theorem}

\begin{theorem}[Proposition 18,~\cite{kantorII}]
Let $(\A;\cdot,\{\cdot,\cdot\})$ be a Poisson algebra.
    Then 
    \[ (\A;\llbracket \cdot+ \{\cdot,\cdot\},\cdot+ \{\cdot,\cdot\} \rrbracket) \]
    is a~noncommutative Jordan algebra.
\end{theorem}

\begin{theorem}[Corollary 31,~\cite{kantorII}]
Let $(\A;\cdot,\circ)$ be a~left Novikov-Poisson algebra.
    Then $(\A,\llbracket \circ,\cdot \rrbracket,\llbracket \cdot,\circ\rrbracket)$ is a~left Novikov-Poisson algebra.
    \end{theorem}

Also, we give constructive methods for obtaining new transposed Poisson algebras
and Poisson-Novikov algebras;
and for
classifying Poisson structures and
commutative post-Lie structures on a~given algebra and construct
some new transposed Poisson algebras.

Let us now consider a~special example of transposed Poisson algebras constructed in~\cite[Theorem 25]{fkl21}.
The transposed Poisson algebra
$(\mathcal{W}, \cdot, [\cdot,\cdot])$ is spanned by the generators
$\{L_i, \ I_j \}_{ i,j \in \mathbb{Z}}$.
These generators satisfy
\begin{longtable}
{rclrcl}
$[L_m, L_n]$&$=$&$ (m -n)L_{m+n} $, & $[L_m, I_n]$&$ =$&$ (m-n - a~)I_{m+n} $,\\
$ L_m \cdot L_n$&$=$&$ w L_{m+n} $, & $L_m \cdot I_n $&$=$&$ w I_{m+n} $,
\end{longtable}
\noindent{}where $w$ is a fixed element from the vector space generated by
$\{L_i, \ I_j \}_{ i,j \in \mathbb{Z}}$ and the multiplication given by juxtaposition satisfies
$L_iL_j=L_{i+j}$ and $L_iI_j=I_{i+j} $.

\begin{theorem}[Proposition 27,~\cite{kantorII}]
Let $\star=\llbracket [\cdot,\cdot],\cdot \rrbracket $ and $\{\cdot,\cdot\}=\llbracket\ \cdot, [\cdot,\cdot] \rrbracket $ be new multiplications defined on multiplications of the transposed Poisson algebra $(\mathcal{W}, \cdot, [\cdot,\cdot])$ defined above.
Then $(\mathcal{W}, \star, \{\cdot,\cdot \})$ is a~transposed Poisson algebra.
\end{theorem}

\subsection{Universal conservative superalgebras}
Kantor introduced the class of conservative algebras in~\cite{Kan1}.
This class includes some well-known classes of algebras, such as associative, Jordan, Lie, Leibniz and Zinbiel~\cite{klp}.
To define conservative algebras we first introduce some notations. Let $V$ be a~vector space over a~field $\mathbb{F} $, let $\varphi$ be a~linear map on $V$, and let $B$ be a~bilinear map on $V$ (i.~e., an algebra). Then we can consider a~product of $\varphi$ and $B$, which is the bilinear map $[\varphi,B]$ on $V$ given by
\begin{equation*}[\varphi,B]
(x,y)= \varphi(B(x,y))- B(\varphi(x),y)-B(x,\varphi(y)).
\end{equation*}
Note that this product measures how far is $\varphi$ from being a~derivation in the algebra $B$.

We say that an algebra $\A$ with a~multiplication $\cdot$ is called a~(left) conservative algebra
if there exists a~(possibly different) algebra structure $* $ (called an associated algebra) on the underlying space of $\A$ such that
\[
[L_a,[L_b, \cdot]]=- [L_{a*b},\cdot].
\]
By replacing the left multiplications with the right multiplications, we can define the right Kantor product and obtain a~similar theory.
In the theory of conservative algebras, the conservative algebra $U(n)$ introduced in Section~\ref{kanpro} has high importance.
One can verify that the algebra $U(n)$ is conservative with the associated multiplication $\ast$ given, for example, by
$A \ast B(x,y) = -B(u,A(x,y))$
(there are other associated multiplications). In the theory of conservative algebras, the algebra $U(n)$ plays a~role analogous to the role of $\mathfrak{gl}_n$ in the theory of Lie algebras, that is, any conservative algebra (modulo its maximal Jacobian ideal) embeds in $U(n)$ for certain $n $.
In 1989, Kantor defined one generalization of conservative algebras,
which he called quasi-conservative algebras
(now known as ``rigid'' algebras, see~\cite{kacan}).
In 2010, Kac and Cantarini considered ``rigid'' (commutative and anticommutative) superalgebras and described simple superalgebras under some special conditions~\cite{kacan}.

One may also use the general approach to define conservative
superalgebras. Namely, let $\Gamma:=\Gamma_0\oplus\Gamma_1$
  be the Grassmann superalgebra in generators $1,\ \xi_i,\ i\in \mathbb N,$
  \[  \Gamma_0=\left<1,\xi_{i_1}\ldots \xi_{i_{2k}}:k\in
  \mathbb N \right>, \Gamma_1=\left<\xi_{i_1}\ldots \xi_{i_{2k-1}}:k\in
  \mathbb N  \right>. \] Let $\A:=\A_0\oplus\A_1 $ be a~superalgebra and $\cdot$ and $*$ be two products on $\A $. Consider its \textit{Grassmann envelope}
 $\Gamma({\A}):=
 (\A_0\otimes\Gamma_0)\oplus(\A_1\otimes\Gamma_1)$ and extend the products $\cdot$ and $*$ to $\Gamma(\A)$ as follows:
 \[
(a\otimes f)\cdot(b\otimes g)=(-1)^{ab}ab\otimes fg,
\]
 \[
(a\otimes f)*(b\otimes g)=(-1)^{ab}a*b\otimes fg
\]
 for all homogeneous $a,b\in {\A},f,g\in \Gamma\ (p(a)=p(f),\,
 p(b)=p(g)) $. Then $(\A,\cdot)$ is conservative with an associated multiplication $*$ if and only if $(\Gamma(\A),\cdot)$ is a~conservative algebra with an associated multiplication $* $.

 The following theorem provides different
 examples of conservative superalgebras.

\begin{theorem}[Theorem 2,~\cite{kpp19}]
Let ${\cal V}$ be a~homogeneous variety of algebras. Assume that
there exist $\alpha, \beta \in \mathbb{F}$
 such that every ${\cal V}$-algebra is conservative with
  the associated multiplication is given by the rule
 $a\,*\,b~=~\alpha ab~+~\beta ba $. Then every ${\cal V}$-superalgebra is
conservative with the associated multiplication $a\,*\,b~=~\alpha
ab + (-1)^{ab}\beta ba $.
\end{theorem}

 It follows that associative, quasi-associative, Jordan, terminal, Lie, Leibniz, and Zinbiel superalgebras are conservative (see~\cite{klp}). In particular, a~superalgebra $M$ is terminal if and only if it is conservative and the multiplication in the associated superalgebra $M^*$ can be given by
\[
M^* (x,y)=\frac{2}{3}xy+(-1)^{xy}\frac{1}{3}yx.
\]

\begin{theorem}[Proposition 4,~\cite{kpp19}]
A conservative superalgebra $\A$ with a~unity is a
noncommutative Jordan superalgebra. \rm
\end{theorem}

Let $V$ be a~superspace. The space of the superalgebra
$U(V)$ is the superspace of all bilinear operations
on $V $.
 Fix a~nonzero homogeneous $a\in V$. Define the multiplication
 ${\Delta} _a$ in $U(V)$ by the rule
\begin{center}
$(A {\Delta} _a B)(x,y) = A(a,B(x,y))- (-1)^{\scriptscriptstyle
 B(A+a)}B(A(a,x),y)- (-1)^{\scriptscriptstyle (A+a)(B+x)}B(x,A(a,y)) $.
\end{center}

Consider the natural action of the group ${\rm GL}(V)$ of
even automorphisms of $V$ on $U(V):$
\begin{center}
$ \varphi(A)(x,y) =
\varphi(A(\varphi^{-1}(x),\varphi^{-1}(y)))$
\end{center}
(note that we denote an automorphism and its action by the same symbol
$\varphi$). A direct computation shows that the mapping $A \mapsto
\varphi(A)$ is an isomorphism between
$(U(V),\Delta_a)$ and
$(U(V),\Delta_{\varphi(a)})$
Therefore, different nonzero
even (respectively, odd) vectors $a$ give rise to isomorphic {\it
even} (respectively, \textit{odd}) superalgebras, which we denote by
 $U(V)^0$ and
 $U(V)^1  $, respectively.
Moreover, consider the opposite superspace $V^{\Pi}$ given by
$V^{\Pi}_{\bar{0}} = V_{\bar{1}}, V^{\Pi}_{\bar{1}} = V_{\bar{0}} $.
Then the parity-reversing isomorphism $V \cong V^{\Pi}$ induces an
 isomorphism between $U(V^{\Pi})^1$ and the odd
superalgebra obtained from $U(V)^0$ by reversing the
parity. Therefore, it suffices to consider only the superalgebras
$U(V)^0  $. For the sake of simplicity, we denote them
  by $U(V) $.
If $V = V_{n,m}$ is a~finite-dimensional superspace with
 ${\rm dim} V_{\bar 0}=n$
and ${\rm dim} V_{\bar 1}=m$ then we denote $U(V)$ by
$U(n,m) $. Further, in this case, we say that $V$ is of
dimension $n+m$.

\begin{theorem}[Theorem 8,~\cite{kpp19}]
Let $V$ be a~superspace, and let $a\in V_{\bar 0}$.
 The superalgebra $(U(V),\Delta_a)$
is conservative, and the associated multiplication can be given by
\begin{center}
$A \bigtriangledown^1_aB(x,y)=-(-1)^{AB}B(a,A(x,y))$
\end{center}
or
\begin{center}
$A \bigtriangledown^2_aB(x,y)= \frac{1}{3}(A^* \Delta_aB + (-1)^{AB}\tilde{B}\Delta_aA) $,
\end{center}
where $A^* (x,y)=A(x,y)+(-1)^{xy}A(y,x)$ and $\tilde{B}(x,y) = 2(-1)^{xy}B(y,x)-B(x,y) $.
\end{theorem}
An even element $e \in \A$ is said to be a~\textit{left quasiunity} if the equality
\[ e(xy)=(ex)y+x(ey)-xy \] holds.

An element $a$ in a~superalgebra $\A$ is called a~\textit{Jacobi
element} provided that
\[ a(xy)=(ax)y+(-1)^{ax}x(ay). \] Let $\A$ be a~conservative superalgebra on a~space $V$ with the
Jacobi subspace $J$. Consider the space $W $, which we define as
$W=V/J$ if $\A$ has a~left quasiunity, and $W=V/J\oplus E$ in the
opposite case, where $E$ is the one-dimensional even space with a
basis element $\epsilon $.

Assume that $\A$ possesses a~quasiunity. Define the \textit{adjoint
mapping} $\operatorname{ad}: \A \to U(W)$ as follows:
\begin{center}
$\operatorname{ad}(a)(\alpha, \beta) =
(-1)^{\beta(a+\alpha)} ((\beta * a) * \alpha +
(-1)^{\alpha a}\beta * (\alpha a)-(-1)^{\alpha a} (\beta * \alpha)* a) $.
\end{center}
 If $\A$ does not have a~quasiunity,
 we define the adjoint mapping
$\operatorname{ad}: \A \to U(W)$ by the equation above and
the following equations:
\begin{center}
$\operatorname{ad}(a)(\alpha, \epsilon)
= a* \alpha + (-1)^{\alpha a}\alpha a~- (-1)^{\alpha a}\alpha * a,
$

$\operatorname{ad}(a)(\epsilon, \beta) =
(-1)^{a\beta}\beta * a,\ \operatorname{ad}(a)(\epsilon, \epsilon)=a $.
\end{center}

\begin{theorem}[Theorem 11,~\cite{kpp19}]
Let $\A$ be a~conservative superalgebra on a~vector space $V$ with the
Jacobi subspace $J $. Let either $W=V/J$ or $W=V/J\oplus
\left<\epsilon\right>$ as above. The adjoint mapping
$\operatorname{ad}: \A \to (U(W),
\Delta_{-\epsilon})$ is a~homomorphism whose kernel is
 the maximal Jacobi ideal. In particular, if $V$ is
finite-dimensional and $J$ is of codimension $n+m $, then we have a
homomorphism $\operatorname{ad}: \A \to U(k,m) $, where $k =
n$ if $\A$ has a~quasiunity and $k = n+1$ otherwise.
\end{theorem}

\subsection{Conservative algebras of $2$-dimensional algebras}
Multiplication on a~2-dimensional vector space is defined by a~$2\times 2\times 2$ matrix. The classification of algebras of dimension $2$ was given in many papers (see, for example,~\cite{kv16}).
Let us consider the space $U(2)$ of all multiplications on the $2$-dimensional space $V_2$ with a~basis $v_1,v_2$.
The definition of the multiplication on the algebra $U(2)$ can be found in Section~\ref{kanpro}. The algebra $U(2)$ is conservative.
Let us consider the multiplications $\alpha_{i,j}^k$ ($i,j,k=1,2$) on $V_2$ defined by the formula $\alpha_{i,j}^k (v_t,v_l)=\delta_{it}\delta_{jl} v_k$ for all $t,l$. It is easy to see that $\{ \alpha_{i,j}^k | i,j,k=1,2 \}$ is a~basis of the algebra $U(2)$.
The multiplication table of $U(2)$ in this basis is given in~\cite{klp}.
But in the next work~\cite{kayvo}, we found another useful basis that will be interesting to consider.
Let's introduce the notation
\begin{longtable}
{llll}
$e_1=\alpha_{11}^1 -\alpha_{12}^2 -\alpha_{21}^2  $, &
$e_2=\alpha_{11}^2  $, &
$e_3=\alpha_{22}^2 -\alpha_{12}^1 -\alpha_{21}^1  $, &
$e_4=\alpha_{22}^1  $,\\
$e_5=2\alpha_{11}^1 +\alpha_{12}^2 +\alpha_{21}^2  $, &
$e_6=2\alpha_{22}^2 +\alpha_{12}^1 +\alpha_{21}^1  $, &
$e_7=\alpha_{12}^1 -\alpha_{21}^1  $, &
$e_8=\alpha_{12}^2 -\alpha_{21}^2  $.
\end{longtable}

It is easy to see that the multiplication table of $U(2)$ in the basis $e_1,\dots,e_8$ is the following.

\begin{center}
\begin{tabular}
{c|c|c|c|c|c|c|c|c|}
      & $e_1$ & $e_2$ & $e_3$ & $e_4$ & $e_5$ & $e_6$ & $e_7$ & $e_8$ \\ \hline
$e_1$ & $-e_1$ & $-3e_2$ & $e_3$ & $3e_4$ & $-e_5$ & $e_6$ & $e_7$ & $-e_8$ \\ \hline
$e_2$ & $3e_2$ & $0$ & $2e_1$ & $e_3$ & $0$ & $-e_5$ & $e_8$ & $0$ \\ \hline
$e_3$ & $-2e_3$ & $-e_1$ & $-3e_4$ & $0$ & $e_6$ & $0$ & $0$ & $-e_7$ \\ \hline
$e_4$ & $0$ & $0$ & $0$ & $0$ & $0$ & $0$ & $0$ & $0$ \\ \hline
$e_5$ & $-2e_1$ & $-3e_2$ & $-e_3$ & $0$ & $-2e_5$ & $-e_6$ & $-e_7$ & $-2e_8$ \\ \hline
$e_6$ & $2e_3$ & $e_1$ & $3e_4$ & $0$ & $-e_6$ & $0$ & $0$ & $e_7$ \\ \hline
$e_7$ & $2e_3$ & $e_1$ & $3e_4$ & $0$ & $-e_6$ & $0$ & $0$ & $e_7$ \\ \hline
$e_8$ & $0$ & $e_2$ & $-e_3$ & $-2e_4$ & $0$ & $-e_6$ & $-e_7$ & $0$ \\ \hline
\end{tabular}
\end{center}

The subalgebra spanned by the elements $e_1, \ldots, e_6$ is the terminal algebra $W_2$ of commutative 2-dimensional algebras.
The subalgebra spanned by the elements $e_1, \ldots, e_4$ is the conservative (and, moreover, terminal) algebra $S_2$ of all commutative 2-dimensional algebras with trace zero multiplication~\cite{klp}.
In a~series of our papers devoted to the study of $U(2)$~\cite{klp, kayvo}, we described
one-sided ideals, subalgebras of codimension $1$, derivations and automorphisms of $U(2) $, $W_2$ and $S_2 $.

We summarize all the main results of the present part in the following theorems.

\begin{theorem}[Theorem 11,~\cite{kayvo}]
${\rm Aut}\big(U(2)\big)\cong {\rm Aut}(W_2)\cong {\rm Aut}(S_2)$. Also, this group is isomorphic to the matrix group $
\begin{pmatrix}
1 & \FF\\ 0 & \FF^*
\end{pmatrix}
$ (the one-dimensional affine group over $\FF$).
\end{theorem}

\begin{theorem}[Theorem 11,~\cite{kayvo}]
The set of nonzero idempotents of the algebra $U(2)$ equals the disjoint union of the following sets (where $c,d\in\FF$, $q\in\bar\FF$):
\begin{longtable}
{l}
 $\mathcal{O}(e_8+e_2-e_1+c(3e_8+e_5-2e_1))$ \\
 $\mathcal{O}(-e_1+c(e_5-2e_1)+de_8)$ \\
 $\mathcal{O}(-e_1-2e_8+4e_3+e_6+3e_7+c(3e_8-e_5+2e_1)+de_4)$ \\
 $\mathcal{O}(-e_1-2e_8+c(3e_8-e_5+2e_1)+qe_4)$.
\end{longtable}
\end{theorem}

\section{Non-associative algebras and superalgebras}
The present part is based on the papers
written together with
Antonio Jesús Calderón,
Artem Lopatin,
Bakhrom Omirov,
Gulkhayo Solijanova,
Luisa Camacho,
Pasha Zusmanovich
and
Yury Popov~\cite{klp18mjm,KZ21,cklos20,ccko21}.

\subsection{Leibniz algebras}
Gradings by abelian groups have played a~key role in the study of Lie algebras and superalgebras, starting with the root space decomposition of the semisimple Lie algebras over the complex field, which is an essential ingredient in the Killing-Cartan classification of these algebras.
Gradings by a~cyclic group appear in the connection between Jordan algebras and Lie algebras through the Tits-Kantor-Koecher construction and in the theory of Kac-Moody Lie algebras.
Gradings by the integers or the integers modulo $2$ are ubiquitous in Geometry.
In 1989, Patera and Zassenhaus~\cite{pz89} began a~systematic study of gradings by abelian groups on Lie algebras. They raised the problem of classifying the fine gradings, up to equivalence, on the simple Lie algebras over the complex
numbers.
After that, gradings of simple alternative and simple Malcev algebras~\cite{elduque98},
the simple Kac Jordan superalgebra~\cite{calderon10},
countless simple Lie algebras~\cite{ draper16, elduque15} and
nilpotent Lie algebras~\cite{bgr16} were described.

In the past years, Leibniz algebras have been under active research (see, for example,~\cite{ leib2,ikv18,kppv, Khudoyberdiyev13, Khudoyberdiyev14}).
Recently, they have appeared in many geometric and physics applications (see, for example,~\cite{bonez, leib2, kotov20,strow20} and references therein).
The main result on the structure of finite-dimensional Leibniz algebras asserts that a~Leibniz algebra decomposes into a~semidirect sum of a~solvable radical and a~semisimple Lie algebra. Therefore, the main problem of the description of finite-dimensional Leibniz algebras consists of the study of solvable Leibniz algebras. Similarly to the case of Lie algebras, the study of solvable Leibniz algebras is reduced to nilpotent ones.
Since the description of all $n$-dimensional nilpotent Leibniz algebras is an unsolvable task (even in the case of Lie algebras), we have to study nilpotent Leibniz algebras under certain conditions (conditions on index of nilpotency, various types of grading, characteristic sequence etc.)~\cite{ Khudoyberdiyev13, Khudoyberdiyev14}. The well-known natural grading of nilpotent Lie and Leibniz algebras is very helpful when investigating the properties of those algebras without restrictions on the grading. Indeed, we can always choose a~homogeneous basis and thus the grading allows us to obtain more explicit conditions for the structural constants. Moreover, such grading is useful for the investigation of cohomologies for the considered algebras, because it induces the corresponding grading of the group of cohomologies. Thus, it is crucial to know what kind of grading a~nilpotent Leibniz algebra admits.

In our work~\cite{ccko21} we begin the study of gradings on Leibniz algebras by classifying, up to equivalence, all abelian groups gradings of null-filiform and one-parametric filiform Leibniz algebras.
Let us define all the considered algebras.

\begin{definition}
An arbitrary complex $n$-dimensional null-filiform Leibniz algebra is isomorphic to the algebra
\[
e_ie_1=e_{i+1}, \quad 1 \leq i \leq n-1.
\]
\end{definition}

\begin{definition}
A complex $n$-dimensional one-parametric filiform Leibniz algebra is isomorphic to the algebra of the following form:
\[
e_1e_2 = \theta e_n, \ e_ie_1 = e_{i+1}, \ 2\leq i \leq {n-1}.
\]
\end{definition}

We described all abelian gradings of null-filiform Leibniz algebras and one-parametric filiform Leibniz algebras in~\cite[Theorems 13, 17 and 19]{ccko21}.

A very powerful tool in the study of nilpotent algebras is the characteristic
sequence, in which a~priori gives the multiplication on one element of the basis.

\begin{definition}
For a~nilpotent Leibniz (Lie) algebra $\A$ and $x\in \A\setminus \A^2$ we consider the decreasing sequence $C(x)=(n_1,n_2,
\ldots,n_k)$ with respect to the lexicographical order of the dimensions Jordan's blocks of the operator $R_x$. The sequence $C(L)=\max\limits_{x\in \A \setminus \A^2}C(x)$
is called the \textit{characteristic sequence} of the Leibniz algebra $\A$.
\end{definition}
Recently, in the paper~\cite{Ancochea} it was considered a~finite-dimensional solvable Lie algebra $\mathfrak{r}_c$ whose nilpotent radical $\mathfrak{n}_c$ has the simplest structure with a~given characteristic sequence $c=(n_1, n_2, \ldots, n_k, 1)$. Using the Hochschild -- Serre factorization theorem, the authors established that the low order cohomology groups of the algebra $\mathfrak{r}_c$ with coefficients in the adjoint representation are trivial.
In our paper~\cite{cklos20}, we study a~family of nilpotent Leibniz algebras whose corresponding Lie algebra is $\mathfrak{n}_c $.
Further, solvable Leibniz algebras with such nilpotent radicals and $(k+1)$-dimensional complementary subspaces to the nilpotent radicals are described. Namely, we prove that such solvable Leibniz algebra is unique and centerless.

\[
\mathbf{R}: \left\{
\begin{array}
{lll}
[e_{1},e_1]=h, \quad [h,x_1]=2h,&\\[1mm]
[e_i,e_1]=-[e_1,e_i]=e_{i+1}, &2 \leq i \leq n_1, &\\[1mm]
[e_{n_1+\ldots+n_{j}+i},e_1]=-[e_1,e_{n_1+\ldots+n_{j}+i}]=e_{n_1+\ldots+n_{j}+1+i},& 2\leq i\leq n_{j+1}, \\[1mm]
[e_1,x_1]=-[x_1,e_1]=e_1, &\\[1mm]
[e_i,x_1]=-[x_1,e_i]=(i-2)e_i,& 3\leq i \leq n_1+1,\\[1mm]
[e_{n_1+\ldots+n_{j}+i},x_1]=-[x_1, e_{n_1+\ldots+n_{j}+i}]=(i-2)e_{n_1+\ldots+n_{j}+i} & 2\le i\leq n_{j+1},\\[1mm]
[e_i,x_{2}]=-[x_{2}, e_i]=e_{i}, & 2\le i\leq n_1+1,\\[1mm]
[e_{n_1+\ldots+n_{j}+i},x_{j+2}]=-[x_{j+2}, e_{n_1+\ldots+n_{j}+i}]=
e_{n_1+\ldots+n_{j}+i}, & 2\le i\leq n_{j+1}.\\[1mm]
\end{array}
\right.
\]
where $1\leq j\leq k-1 $.
For this Leibniz algebra, the triviality of the first and the second cohomology groups with coefficients in the adjoint representation is established as well.

\begin{theorem}[Theorem 21,~\cite{cklos20}]
The first and second groups of cohomologies
of the solvable Leibniz algebra $\mathbf{R}$ in coefficient itself are trivial.
\end{theorem}

\begin{theorem}[Corollary 22,~\cite{cklos20}]
 The solvable Leibniz algebra $\mathbf{R}$ is rigid.
\end{theorem}

\begin{remark}
Note that the structure of the rigid algebra $\mathbf{R}$ depends on the given decreasing sequence $(n_1, n_2, \ldots, n_k) $. Set $p(n)$ to be the number of such sequences, that is, $p(x)$ is
the number of integer solutions of the equation
\begin{center}
    $n_1 + n_2 + \ldots + n_k = n$ with $n_1 \geq n_2 \geq \ldots \geq n_k \geq 0 $.
\end{center}
The asymptotic value of $p(n)$, given by the expression
\begin{center}
    $p(n)\approx\frac{1}{4n\sqrt{3}} e^{\pi\sqrt{2n/3}}$ \ (where $a(n)\approx b(n)$
means that $\lim\limits_{n\rightarrow\infty}\frac{a(n)}{b(n)}=1$),
\end{center}
yields the existence of at least $p(n)$ irreducible components of the variety of Leibniz algebras of dimension $n+k+3$.

\end{remark}

\subsection{Anticommutative $\mathfrak{CD}$-algebras}

The idea of the definition of the variety of $\mathfrak{CD}$-algebras comes from the following property of Jordan and Lie algebras: \textit{the commutator of any pair of multiplication operators is a~derivation}.
 This gives three identities of degree four, which reduces to only one identity of degree four in the commutative or anticommutative case.
Commutative and anticommutative $\mathfrak{CD}$-algebras are related to many interesting varieties of algebras.
 Thus, anticommutative $\mathfrak{CD}$-algebras are a~generalization of Lie algebras,
 containing the intersection of Malcev and Sagle algebras as a~proper subvariety. Moreover, the following intersections of varieties coincide:
Malcev and Sagle algebras;
Malcev and anticommutative $\mathfrak{CD}$-algebras; and
Sagle and anticommutative $\mathfrak{CD}$-algebras.
On the other hand,
the variety of anticommutative $\mathfrak{CD}$-algebras is a~proper subvariety of
the varieties of binary Lie algebras
and almost Lie algebras~\cite{KZ21}.
The variety of anticommutative $\mathfrak{CD}$-algebras coincides with the intersection of the varieties of binary Lie algebras and almost Lie algebras.
Commutative $\mathfrak{CD}$-algebras are a~generalization of Jordan algebras,
which are a~generalization of commutative associative algebras.
On the other hand, the variety of commutative $\mathfrak{CD}$-algebras is also known as the variety of almost-Jordan algebras, which lies in the bigger variety of generalized almost-Jordan algebras~\cite{fl15}.
 The $n$-ary version of commutative $\mathfrak{CD}$-algebras was introduced in a~recent paper by
Kaygorodov, Pozhidaev and Saraiva~\cite{KPS19}.
The variety of almost-Jordan algebras is the variety of commutative algebras
satisfying
\[
2((yx)x)x+yx^3 =3(yx^2)x.
\]
This present identity appeared in 1965 in a~paper by Osborn
during the study of identities of degrees less than or equal to $4$ of non-associative algebras. The identity is a~linearized form of the Jordan identity.
The systematic study of almost-Jordan algebras was initiated in the following paper by Osborn and it was continued in some papers by Petersson~\cite{petersson67}, Osborn~\cite{osborn69}, and Sidorov
(sometimes, they were called Lie triple algebras).
Hentzel and Peresi proved that every semiprime almost-Jordan ring is Jordan~\cite{peresi}.
After that,
Labra and Correa
proved that a~finite-dimensional almost-Jordan right-nilalgebra is nilpotent.
Assosymmetric algebras under the symmetric product give almost-Jordan algebras~\cite{askar18}.

An anticommutative algebra $\A$ is called a~$\mathfrak{CD}$-algebra if it satisfies the property that
for any $a,b\in \A$, the commutator $[R_a, R_b]$ is a~derivation of $\A$. This
condition can be written as a~homogeneous identity of degree $4$ comprising $6$
monomials:
\[
((xy)a)b - ((xy)b)a - ((xa)b)y + ((xb)a)y + ((ya)b)x - ((yb)a)x = 0.
\]
We introduced the notion of anticommutative $\mathfrak{CD}$-algebras and studied their properties in~\cite{KZ21}.
Firstly, we compare the variety of anticommutative $\mathfrak{CD}$-algebras with other known varieties, namely
Lie, binary Lie, Malcev, and Sagle.
Secondly, we consider the cohomology theory of anticommutative $\mathfrak{CD}$-algebras and obtain some relations with Lie algebras. Namely, we prove the following results.

\begin{theorem}[Proposition 3,~\cite{KZ21}]
For any $\mathfrak{CD}$-algebra $\A$, the quotient $\A/Z(\A)$ is a~Lie algebra. In particular,
any centerless (and, in particular, simple) $\mathfrak{CD}$-algebra is a~Lie algebra.
\end{theorem}


\begin{theorem}[Theorem,~\cite{KZ21}]
For any simple finite-dimensional Lie algebra $\A$ over a~field of characteristic
zero, and any finite-dimensional $\A$-module $M$, ${\rm H}_{\mathfrak{CD}}^2 (\A,M) = 0$.
\end{theorem}

\subsection{Noncommutative Jordan superalgebras}
The class of noncommutative Jordan algebras is vast. 
As an example, it includes alternative algebras, Jordan algebras, quasiassociative algebras, quadratic flexible algebras, and anticommutative algebras.
Schafer proved that a~simple noncommutative Jordan algebra is either a~simple Jordan algebra, a~simple quasiassociative algebra, or a~simple flexible algebra of degree $2$.
Oehmke proved the analog of Schafer's classification for simple flexible algebras with strictly associative powers in the case of characteristic different from $2$, $3$,
McCrimmon classified simple noncommutative Jordan algebras of degree greater than $2$ in the case of characteristic different from $2$, and
Smith described such algebras of degree $2$.
The case of nodal simple algebras of positive characteristic was considered in the papers by Kokoris.
The case of simple finite-dimensional Jordan superalgebras over algebraically closed fields of characteristic zero was studied by Kac and Kantor.
Racine and Zelmanov classified the finite-dimensional Jordan superalgebras of characteristic different from $2$ with a~semisimple even part. The case when the even part is not semisimple was considered by Mart\'{i}nez and Zelmanov in~\cite{MZ01}
and Cantarini and Kac described all linearly compact simple Jordan superalgebras~\cite{KacKant07}.
Simple noncommutative Jordan superalgebras were described by Pozhidaev and Shestakov in~\cite{PS13,PS19}.
Representations of simple noncommutative superalgebras were described by Popov~\cite{Po1}.
Nowadays, the study of properties of simple non-associative algebras and superalgebras has aroused a~strong interest.
For example, Popov determined the structure of differentiably simple Jordan algebras, and
Barreiro, Elduque, and Mart\'{i}nez described the derivations of the Cheng-Kac Jordan superalgebra~\cite{BEC}. Moreover,
Kaygorodov, Shestakov, and Zhelyabin studied generalized derivations of Jordan algebras and superalgebras~\cite{lesha12,lesha14}.
Another interesting problem in the study of Jordan algebras and superalgebras is a~description of maximal subalgebras and automorphisms~\cite{ELS,ELS2}.

A superalgebra $\A$ is called a~noncommutative Jordan superalgebra if it satisfies the following operator identities:
\[
[R_{x \circ y}, L_z] + (-1)^{|x|(|y|+|z|)}[R_{y\circ z}, L_x] + (-1)^{|z|(|x|+|y|)}[R_{z \circ x}, L_y]=0,
\]
\[
[R_x, L_y] = [L_x, R_y].\label{flex}
\]

A binary linear operation $\{ \cdot,\cdot\}$ is called a generic Poisson bracket~\cite{ksu18} on a~superalgebra $(\A,\cdot)$ if for arbitrary homogeneous $a, b, c \in \A$ we have
\[
\{a \cdot b, c\} = (-1)^{|b||c|}\{a,c\}\cdot b + a~\cdot \{b,c\}.
\]
We notice that there is a~one-to-one correspondence between noncommutative Jordan superalgebras and superanticommutative Poisson brackets on adjoint Jordan superalgebras.
Let $(\A, \bullet)$ be a~Jordan superalgebra and $[\cdot,\cdot]$ be a generic Poisson bracket on $\A $. Then the operation $ab = \frac{1}{2}(a\bullet b + [a,b])$ turns $(\A, \bullet,[\cdot,\cdot])$ into a~noncommutative Jordan superalgebra. Conversely, if $\A$ is a~noncommutative Jordan superalgebra, then the supercommutator $[\cdot,\cdot]$ is a~generic Poisson bracket on a~Jordan superalgebra $\A^{(+)} $. Moreover, the multiplication in $\A$ can be recovered by the Jordan multiplication in $\A^{(+)}$ and the generic Poisson bracket $[\cdot,\cdot]: ab = \frac{1}{2}(a \bullet b + [a,b]) $.

In our paper with Lopatin and Popov, we study simple noncommutative Jordan superalgebras constructed in some
papers by Pozhidaev and Shestakov~\cite{PS13}.
We describe all
subalgebras and automorphisms of simple noncommutative Jordan superalgebras $K_3(\alpha,\beta,\gamma)$ and $D_t(\alpha,\beta,\gamma)$~\cite{klp18mjm}. We also
compute the derivations of the nontrivial simple finite-dimensional noncommutative Jordan superalgebras~\cite[Theorems 9 and 10]{klp18mjm}].

\section{Maps on associative and non-associative algebras}
The present part is based on the papers
written together with
Artem Lopatin,
Bruno Ferreira,
Feng Wei,
Henrique Guzzo,
Mykola Khrypchenko
and
Yury Popov~\cite{KLP18inv,KKW19,FK21alt,FGK21alt,KP}.

\subsection{Alternative and Jordan algebras with invertible derivations}

In 1983, Bergen, Herstein, and Lanski initiated a~study whose purpose was to relate the structure of a~ring to the special behavior
of one of its derivations. Namely, in their article~\cite{Berg} they described associative rings admitting derivations with invertible values,
i.e., there is a~derivation $d$ such that $d(a)$ is an invertible element or zero for every $a$ from our algebra.
They proved that such ring must be either a~division ring, or the ring of $2 \times 2$ matrices over a~division ring,
or a~factor of a~polynomial ring over a~division ring of characteristic 2. They also characterized the division rings
such that the $2 \times 2$ matrix ring over them has an inner derivation with invertible values.
Further, associative rings with derivations with invertible values (and also their generalizations) were discussed in a~variety of works
(see, for instance,~\cite{Giam}).
So, in~\cite{Giam} semiprime associative rings with involution, allowing a~derivation with invertible values
on the set of symmetric elements, were examined.
In their work, Bergen and Carini studied associative rings, admitting a~derivation with invertible values on some non-central Lie ideal.
Also, in the papers by Chang, Hongan and Komatsu, the structure of associative rings that admit $\alpha$-derivations with invertible values and their natural generalizations --- $(\sigma,\tau)$-derivations with invertible values, was described.
Komatsu and Nakajima described associative rings that allow generalized derivations with invertible values.
The case of associative superalgebras with derivations with invertible values was studied in the paper by Demir, Alba\c{s}, Arga\c{c}, and Fosner.

The description of non-associative algebras admitting derivations with invertible values began in our paper with Popov~\cite{KP}, where it was proved that every alternative (non-associative) algebra admitting derivation with invertible values
is a~Cayley --- Dickson division algebra over their center or a~factor-algebra of a~polynomial algebra $C[x]/(x^2)$ over a~Cayley --- Dickson division algebra.
Our research was continued in a~joint paper with Lopatin and Popov~\cite{KLP18inv},
where we describe Jordan algebras admitting derivations with invertible values.
Let us summarize the main results from~\cite{KLP18inv}.

\begin{theorem}[Theorem 3,~\cite{KLP18inv}]
Let $\A$ be a~Jordan algebra admitting a~derivation with invertible values $d$. Then $\A$ is an extension of a~simple Jordan algebra with derivation with invertible values by an ideal $M$ such that $M^2 = 0, d(M) = 0$ and $M$ is the largest ideal of $\A $.
\end{theorem}

\begin{theorem}[Theorem 11,~\cite{KLP18inv}]
Let $\A$ be a~Jordan algebra of characteristic not $2$ admitting a~derivation with invertible values $d$.
Then one of the following holds:
\begin{itemize}
\item[$(1)$] $\A$ is an algebra $A^{(+)}$, where $A = D$ or $D_2$, and $D$ is an associative division algebra;
\item[$(2)$] $\A$ is an algebra $H(A, *)$, where $A = D$ or $D_2$, and $D$ is an associative division algebra;
\item[$(3)$] $\A$ is an algebra of symmetric nondegenerate bilinear form $J(V,f)$;
\item[$(4)$] $\A$ is a~division algebra of Albert type;
\item[$(5)$] $\A$ is an extension of cases $(1)$ -- $(4)$ by $M = \mathbf{P}(\A)$ (where $\mathbf{P}(\A)$ stands for the prime radical of $\A$), where $M \subseteq \ker d$ and $M$ is the largest ideal of $\A$.
\end{itemize}

\end{theorem}

\begin{theorem}[Theorem 20,~\cite{KLP18inv}]
Let $\A$ be a~finite-dimensional Jordan algebra admitting a~derivation with invertible values. Then $\A$ is either simple or an algebra of a~symmetric bilinear form (possibly degenerated).
\end{theorem}

\subsection{Maps on alternative algebras}
The study of Lie isomorphisms of rings was originally inspired by Herstein's generalization of generalizing classical theorems on the Lie structure of total matrix rings to the Lie structure of arbitrary simple rings~\cite{iher}.
In this paper, Herstein formulated some open questions that gave rise to big research in Lie and Jordan maps on various structures.
The main questions from Herstein's program were solved for associative algebras in papers by
Beidar, Brešar, Chebotar and Martindale~\cite{bbcm1,bbcm2}.
In~\cite{Mart}, Martindale studied Lie isomorphisms between primitive rings $\A$, $\A'$, where he assumed that the characteristic of $\A$ is different from $2$ and $3$ and that $\A$ contains three nonzero orthogonal idempotents whose sum is the identity. A few years later, he studied Lie isomorphisms between two simple rings $\A$, $\A'$.
Alternative algebras, as well as associative algebras, admit a~Peirce decomposition.
The recently introduced axial algebras reflect an active interest in the study of algebras admitting a~Peirce-type decomposition~\cite{kms}.
In our works together with Ferreira and Guzzo,
we continue the study of the Herstein program for some special types of alternative algebras admitting Peirce decomposition.

In the following theorems, we summarize the main results of the present part.

\begin{theorem}[Theorem 2.1,~\cite{FGK21alt}]
Let $\A$ be a~unital prime (simple, simple associative) alternative algebra,
let $\A'$ be another prime (simple, simple associative) alternative algebra and
let $\varphi: \A \rightarrow \A'$ be a~surjective Lie multiplicative map that preserves idempotents. Assume that $\A$ has a~nontrivial idempotent $e_1$ with associated Peirce decomposition
$\A = \A_{11} \oplus \A_{12} \oplus \A_{21} \oplus \A_{22}$, such that
\begin{itemize}
\item [\it (1)] if $x_{ij}\A_{ji} = 0$, then $x_{ij} = 0$ $(i \neq j) $;
\item [\it (2)] if $x_{11}\A_{12} = 0$ or $\A_{21}x_{11} = 0$, then $x_{11} = 0 $;
\item [\it (3)] if $\A_{12}x_{22} = 0$ or $x_{22}\A_{21} = 0$, then $x_{22} = 0 $;
\item [\it (4)] if $z \in Z(\A)$ with $z \neq 0$, then $z\A = \A$.
\end{itemize}
Then, the following holds.
\begin{itemize}
    \item If $f_i\varphi(\A_{jj})f_i \subseteq Z(\A') f_i $,
then $\varphi$ is of the form $\psi + \tau$, where $\psi$ is an additive isomorphism between $\A$ and $\A'$ and $\tau$ is a~map from $\A$ to $Z(\A')$, which maps commutators into zero.

\item If $f_i\varphi(\A_{ii})f_i \subseteq Z(\A') f_i $,
then $\varphi$ is of the form $\psi + \tau$, where -$\psi$ is additive anti-isomorphism between $\A$ and $\A'$ and $\tau$ is a~map from $\A$ to $Z(\A')$, which maps commutators into zero. Observe
$f_i = \varphi(e_i)$ and $f_j = 1_{\A'} - f_i$, $i \neq j$.
\end{itemize}

\end{theorem}

\begin{definition}
Let be $\varphi$ an additive map on $\A$. We call $\varphi$ a commuting map if
$[\varphi(x),x] = 0$ for all $x \in \A$.

\end{definition}

\begin{theorem}[Theorem 4,~\cite{FK21alt}]
Let $\A$ be a~unital alternative algebra. Assume that $\A$ has a~nontrivial idempotent $e_1$ with associated Peirce decomposition $\A = \A_{11} \oplus \A_{12} \oplus \A_{21} \oplus \A_{22}$, such that $(x \A) e_i = {0}$ implies $x = 0$ $\left(i = 1,2\right)$.
Let $\varphi: \A \rightarrow \A$ be an additive map. Then the following statements are equivalent:
\begin{itemize}
\item[$(\spadesuit)$] $\varphi$ is commuting;
\item[$(\clubsuit)$] there exists $z \in Z(\A)$ and an additive map
$\Xi : \A \rightarrow Z(\A)$ such that $\varphi(x) = zx + \Xi(x)$ for all $x \in \A$.
\end{itemize}
\end{theorem}

\subsection{Higher derivations of incidence algebras}\label{secincidence}
	Let $(P,\le)$ be a~preordered set and $R$ be a~commutative unital ring.
	Assume that $P$ is locally finite, i.e. for any $x\le y$ in $P$ there are only finitely many $z\in P$ such that $x\le z\le y$.
	The incidence algebra $I(P,R)$ of $P$ over $R$ is the set of functions
	
\[
	\{f: P\times P\to R\mid f(x,y)=0\ \text{if}\ x\not\le y\}
\]
	with the natural structure of an $R$-module and multiplication given by the convolution
	
\begin{align*}
	(fg)(x,y)&=\sum_{x\le z\le y}f(x,z)g(z,y)	
\end{align*}
	for all $f,g\in I(P,R)$ and $x,y\in P$. It would be helpful to point out that the full matrix algebra $M_n(R)$, as well as the upper triangular matrix algebra $T_n(R)$, are particular cases of incidence algebras. In addition, in the theory of operator algebras, the incidence algebra $I(P,R)$ of a~finite poset $P$ is referred to as a~bigraph algebra or a~finite-dimensional commutative subspace lattice algebra.
	Incidence algebras appeared in the early work by Ward~\cite{Ward} as generalized algebras of arithmetic functions. Later, they were extensively used as the fundamental tool of enumerative combinatorics in the series of works ``On the foundations of combinatorial theory''. The study of algebraic mappings on incidence algebras was initiated by Stanley~\cite{Stanley70}. Since then, automorphisms, involutions, and derivations (and their generalizations) on incidence algebras have been actively investigated, see~\cite{inc1,inc2,inc3} and the references therein.

	Let $P$ be a~poset and $R$ a~commutative unital ring. Recall from~\cite{KN09} that a finitary series is a~formal sum of the form
	
\[
\alpha=\sum_{x\le y}\af_{xy}e_{xy},
\]
	where $x,y\in P$, $\af_{xy}\in R$ and $e_{xy}$ is a~symbol, such that for any pair $x<y$ there exists only a~finite number of $x\le u<v\le y$ with $\af_{uv}\ne 0$. The set of finitary series, denoted by $FI(P,R)$, possesses a~natural structure of an $R$-module. Moreover, it is closed under the convolution
	
\[
	\alpha\beta=\sum_{x\le y}\left(\sum_{x\le z\le y}\alpha_{xz}\beta_{zy}\right)e_{xy}.
\]
	Thus, $FI(P,R)$ is an algebra, called the finitary incidence algebra of $P$ over $R$. The identity element of $FI(P,R)$ is the series $\delta=\sum_{x\in P}1_Re_{xx}$. Here and in what follows we adopt the following convention: if in the formal sum the indices run through a~subset $X$ of the ordered pairs $(x,y)$, $x,y\in P$, $x\le y$, then $\alpha_{xy}$ is meant to be zero for $(x,y)\not\in X$.

Higher
	derivations are closely related to derivations. It should be remarked that the first component $d_1$ of each higher derivation $D=\{d_n\}_{n=0}^\infty$ of an algebra $\A$ is itself a~derivation of $\A$. Conversely, let $d: \A\to\A$ be an ordinary derivation of an algebra $\A$ over a~field of characteristic zero. Then, $D=\{\frac{1}{n!}d^n \}_{n=0}^\infty$ is a~higher derivation of $\A$. Heerema, Mirzavaziri and Saymeh independently proved that each higher derivation of an algebra $\A$ over a~field of characteristic zero is a~combination of compositions of derivations, and hence one can characterize all higher derivations on $\A$ in terms of the derivations on $\A$. Ribenboim systemically studied higher derivations of arbitrary rings and those of arbitrary modules, where some familiar properties of derivations are generalized to the case of higher derivations. Ferrero and Haetinger found the conditions under which Jordan higher derivations (or Jordan triple higher derivations) of a~$2$-torsion-free (semi-)prime ring are higher derivations, and the same authors studied higher derivations on (semi-)prime rings satisfying linear relations. Wei and Xiao~\cite{Wei-Xiao11} described higher derivations of triangular algebras and related mappings, such as inner higher derivations, Jordan higher derivations, Jordan triple higher derivations, and their generalizations.

\begin{definition}
A sequence $d=\hd d$ of additive maps $R\to R$ is a~\emph{higher derivation} of $R$ \emph{(of infinite order)} if it satisfies the conditions
	
\begin{center}
		1. $d_0={\rm id}_R$; \ \ \ \
		2. $d_n(rs)=\sum_{i+j=n}d_i(r)d_j(s)$\label{d_n(rs)=sum-d_i(r)d_j(s)}
	
\end{center}
	for all $n\in\NN=\{1,2,\dots\}$ and all $r,s\in R$. If \cref{d_n(rs)=sum-d_i(r)d_j(s)} holds for all $1\le n\le N$, then the sequence $\{d_n\}_{n=0}^N$ is called a~ \emph{higher derivation of order $N$}. Evidently, $\hd d$ is a~higher derivation if and only if $\{d_n\}_{n=0}^N$ is a~higher derivation of order $N$ for all $N\in\NN$. In particular, $d_1$ is always a~usual derivation of $R$.
\end{definition}

Denote by $\hder R$ the set of higher derivations of $R$ and consider the following operation on $\hder R$
\[
	(d'*d'')_n=\sum_{i+j=n}d'_i\circ d''_j.
\]
In particular,
\[
	(d'*d'')_1=d'_1+d''_1.
\]
It was proved in~\cite{Heerema68} that $\hder R$ forms a~group with respect to $*$, whose identity is the sequence $\hd \e$ with $\e_0=id_R$ and
$ 	\e_n=0$ for $n\in\NN$.

Given $r\in R$ and $k\in\NN$, define
\begin{align*}[r,k]
_0&=id_R,\\
	[r,k]_n(x)&=
\begin{cases}
	 0, & k\nmid n,\\
	 r^l x-r^{l-1}xr, & n=kl,
\end{cases}
\end{align*}
for all $n\in\NN$ and $x\in R$. It was proved in~\cite{Nowicki84-ider} that $\hd{[r,k]}\in\hder R$, so that for any sequence $\{r_n\}_{n=1}^\infty\sst R$ one may define $\hd{(\Dt_r)}$ by means of $(\Dt_r)_0=id_R$ and
	$(\Dt_r)_n=([r_1,1]*\dots *[r_n,n])_n $,
where $n\in\NN$. Higher derivations of the form $\Dt_r$ are called \emph{inner}. By~\cite[Corollary 3.3]{Nowicki84-ider} the set of inner higher derivations forms a~normal subgroup in $\hder R$, which will be denoted by $\ihder R$. In particular,
\[
	(\Dt_r)_1(x)=[r_1,1]_1(x)=r_1x-xr_1
\]
is the usual inner derivation of $R$ associated with $r_1\in R$, which we denote by $\ad_{r_1}$.

\begin{definition}
		A sequence $\s=\hd{\s}$ of maps on $I=\{(x,y)\in P\times P\mid x\le y\}$ with values in $R$ is called a~\emph{higher transitive map} if
		
\begin{itemize}
			\item $\s_0(x,y)=1_R$ for all $x\le y$;\label{s_0(x_y)=1-for-all-x_y}
			\item $\s_n(x,y)=\sum_{i+j=n}\s_i(x,z)\s_j(z,y)$ for all $x\le z\le y$.\label{s_n(x_y)=sum-s_i(x_z)s_j(z_y)}
		
\end{itemize}
	
\end{definition}

\begin{theorem}[Lemma 2.6,~\cite{KKW19}]
		Given a~higher transitive map $\s$, denote by $\tl\s$ the following sequence of maps $FI(P,R)\to FI(P,R)$:
		
\begin{align*}
		\tl\s_n(\af)=\sum_{x\le y}\s_n(x,y)\af_{xy}e_{xy},	
\end{align*}
		where $n\in\NN\cup\{0\}$ and $\af\in FI(P,R)$. Then, $\tl\s\in\hder{FI(P,R)}$.
	
\end{theorem}

\begin{theorem}[Theorem 2.8,~\cite{KKW19}]
		Every $R$-linear higher derivation of $FI(P,R)$ is of the form $\Dt_\rho*\tl\s$ for some $\rho=\{\rho_n\}_{n=1}^\infty\subseteq FI(P,R)$ and some higher transitive map $\s$.
	
\end{theorem}

\section{Generalized derivations of non-associative algebras}

The present part is based on the papers
written together with
Bruno Ferreira,
Inomjon Yuldashev,
Karimbergan Kudaybergenov,
Patrícia Beites,
and Yury Popov~\cite{KKY,kay14mz, kay14st,bkp19,kaypopov16,KP, KP16,FKK}.

\subsection{Local derivations of $n$-ary algebras}
The study of local derivations
started in 1990 with Kadison’s article~\cite{k90}. A similar notion, which characterizes non-linear generalizations of derivations,
was introduced by \v{S}emrl as $2$-local derivations.
In his paper~\cite{Semrl97}, it was
proved that a~$2$-local derivation of the algebra
$B(H)$ of all bounded linear operators on the infinite-dimensional
separable Hilbert space $H$ is a~derivation.
After these works, numerous new results related to the description of local and $2$-local derivations of associative algebras appeared (see, for example,~\cite{Khrypchenko19}).
The study of local and $2$-local derivations of non-associative algebras was initiated in some
papers by Ayupov and Kudaybergenov (for the case of Lie algebras, see~\cite{ak17,ak16}).
In particular, they proved that there are no pure local and $2$-local derivations on semisimple finite-dimensional Lie algebras.
Also,~\cite{AyuKudRak16} provides examples of $2$-local derivations on nilpotent Lie algebras that are not derivations.
After the cited works, the study of local and $2$-local derivations was continued for Leibniz algebras and Jordan algebras~\cite{aa17}.
Local automorphisms and $2$-local automorphisms were also studied in many cases;
for example, they were studied on Lie algebras~\cite{ak17,c19}.

The description of local and $2$-local derivations of $n$-ary algebras started in our paper with Ferreira and Kudaybergenov~\cite{FKK}.
Namely, we gave the first example of a~complex simple finite-dimensional (ternary) algebra with non-trivial local derivations.
The idea of introducing a~generalization of Filippov algebras comes from binary Malcev algebras and it was carried out in a~paper by Pozhidaev~\cite{app}.
He defined $n$-ary Malcev algebras, generalizing Malcev algebras and $n$-ary Filippov algebras.
Let us summarize the construction of the most important example of $n$-ary Malcev (non-Filippov). We denote by $\A$ a~composition $8$-dimensional algebra with an involution $\bar{} : a~\mapsto \bar{a}$ and unity $1$.
The symmetric bilinear form $\left\langle x,y\right\rangle = \frac{1}{2} (x\bar{y} + y\bar{x})$ defined on $\A$ is assumed to be
nonsingular.
If
$\A$ is equipped with a~ternary multiplication $[\cdot,\cdot,\cdot]$ by the rule
\begin{center}
$[x, y, z] = (x\bar{y})z - \left\langle y, z\right\rangle x + \left\langle x, z \right\rangle y - \left\langle x, y \right\rangle z $,
\end{center}
then $\A$ becomes a~ternary Malcev algebra~\cite{app},
which will be denoted by ${\rm M}_8$.

Let us give two main definitions necessary for this part.

\begin{definition}
Let $\A$ be an $n$-ary algebra.
A linear map $\nabla : \A \rightarrow \A$ is called a~local derivation if for
any element $x \in \A$ there exists a~derivation ${\mathfrak D}_x: \A \rightarrow \A$ such that $\nabla(x) = {\mathfrak D}_x(x)$.
\end{definition}
\begin{definition}
A (not necessary linear) map $\Delta: \A \rightarrow \A$ is called a~$2$-local derivation, if for
any two elements $x, y \in \A$ there exists an derivation ${\mathfrak D}_{x,y} : \A \rightarrow \A$ such that
$\Delta(x) = {\mathfrak D}_{x,y}(x)$, $\Delta(y) = {\mathfrak D}_{x,y}(y)$.
\end{definition}

  In the following theorems, we summarize all the main results of the present part.

\begin{theorem}[Theorem 3,~\cite{KKY}]
Every local derivation of a~complex finite-dimensional semisimple Leibniz algebra is a~derivation.
\end{theorem}

\begin{theorem}[Theorems 4 and 5,~\cite{FKK}]
Every local (and $2$-local) derivation of a~complex finite-dimensional simple $n$-Lie algebra is a~derivation.
\end{theorem}

\begin{theorem}[Theorem 8,~\cite{FKK}]
A linear mapping $\nabla$ on $M_{8}$ is a~local derivation if and only if its matrix is antisymmetric.
In particular, the dimension of the space
$LocDer M_{8}$ of all local derivations of ${\rm M}_8$ is equal to $28 $.
\end{theorem}

\subsection{Leibniz-derivations of non-associative algebras}

The theory of Lie algebras having a~nonsingular derivation has a~rich history and is still an active research area.
Such Lie algebras appear in many different situations, such as in the studies of pro-$p$ groups of finite coclass by
Shalev and Zelmanov~\cite{Shalev94,SZ92} and in the problems concerning the existence of left-invariant affine structures on Lie groups (see Burde's survey~\cite{Burde06} for details).

In 1955, Jacobson~\cite{Jac} proved that a~finite-dimensional Lie algebra over a
field of characteristic zero admitting a~nonsingular (invertible) derivation is nilpotent. The problem of whether the inverse of this
statement is correct remained open until the work~\cite{Dix}, where
an example of a~nilpotent Lie algebra in which all derivations are
nilpotent (and hence, singular) was constructed.
For Lie algebras in prime characteristic, the situation is more complicated. In that
case, there exist non-nilpotent Lie algebras, even simple ones, which admit nonsingular derivations~\cite{BKK95}.
The main examples of nonsingular derivations are periodic derivations.
Kostrikin and Kuznetsov~\cite{KK} noted that a~Lie algebra admitting a~nondegenerate derivation admits a~periodic derivation, that is, a~derivation $d$ such that $d^N = {\rm id}$ for some $N$, and proved that a~Lie algebra admitting a~derivation of period $N$ is abelian provided that $N \not \equiv 0 \text{ (mod 6)}$.
Burde and Moens proved that a~finite-dimensional complex Lie algebra $\A$ admits a~periodic derivation if and only if $\A$ admits a~nonsingular derivation whose inverse is again a~derivation if and only if $\A$ is hexagonally graded~\cite{BurdeMoens12}.
In the case of positive characteristic $p$, Shalev proved that if a~Lie algebra $\A$
admits a~nonsingular derivation of order $n=p^s m $, where $(m,p)=1$ and $m<p^2 -1 $,
then $\A$ is nilpotent.
The study of periodic derivations was continued by Mattarei~\cite{mat02,mat07,mat09}.

The study of derivations of Lie algebras led to the appearance of the notion of their
natural generalization --- a~pre-derivation of a~Lie algebra, which is a~derivation of a~Lie triple system induced by that algebra.
In~\cite{Baj} it was proved that Jacobson's result remains true in terms of pre-derivations.
Several examples of nilpotent Lie
algebras whose pre-derivations are nilpotent were presented in~\cite{Baj}.

 A generalization of derivations and
pre-derivations of Lie algebras is defined in the paper~\cite{Moens} as a~Leibniz-derivation of
order $k$.
Moens proved that a~finite-dimensional Lie algebra over a~field of characteristic zero is nilpotent if and only if
it admits an invertible Leibniz-derivation.
After that,
Fialowski, Khudoyberdiyev and Omirov~\cite{Fialc12} showed
that with the definition of Leibniz-derivations from~\cite{Moens},
the analogous result for non-Lie Leibniz algebras is not true.
Namely, they gave an example of non-nilpotent Leibniz algebra, which
admits an invertible Leibniz-derivation. To extend the
results of the paper~\cite{Moens} for Leibniz algebras they introduced the
definition of Leibniz-derivation of a~Leibniz algebra which is coherent
with the definition of Leibniz-derivation of a~Lie algebra and proved
that a~finite-dimensional Leibniz algebra is nilpotent if and only if it admits an
invertible Leibniz-derivation.
In the paper~\cite{KP}, the authors showed that the same result holds for alternative algebras (particularly, for associative algebras).
Also, in this article an example of a~nilpotent alternative (non-associative) algebra over a~field of positive characteristic possessing only singular derivations was provided.

It is well known that the radicals of finite-dimensional algebras belonging to the classical varieties (such as varieties of Jordan algebras, Lie algebras, alternative algebras, and many others) are invariant under their derivations~\cite{Slinko1972}. Therefore, it is natural to state another interesting problem: is the radical of an algebra invariant under its Leibniz-derivations? Moens proved that the solvable radical of a~Lie algebra is invariant under all its Leibniz-derivations, Fialowski, Khudoyberdiev and Omirov showed the invariance of solvable and nilpotent radicals of a~Leibniz algebra, and the authors of~\cite{KP} proved analogous results for alternative algebras.
Another interesting task is to describe the Leibniz-derivations of algebras belonging to certain ``nice'' classes, such as semisimple and perfect algebras. Moens described all Leibniz-derivations of semisimple Lie algebras~\cite{Moens}
and Zhou described all pre-derivations of perfect centerless Lie algebras of characteristic $\neq 2$.

\begin{definition}
Let $\A$ be an algebra, $n$ be a~natural number $ \geq 2$, and $f$ be an arrangement of brackets on a~product of length $n$. A linear mapping $d$ on $\A$ is called an $f$-Leibniz derivation of $\A$, if for any $a_1, \dots a_n \in \A$ we have
\begin{equation*}
 d([a_1,\ldots,a_n]_f) = \sum\limits_{i=1}^n [a_1,a_2,\ldots, d(a_i), \ldots a_n]_f.
\end{equation*}
\end{definition}

Particularly, if $f = l(n)$ ($f = r(n)$) is the left (right) arrangement of brackets of length $n$, that is,
\[
[x_1, \ldots, x_n]_{l(n)} = ((\ldots(x_1x_2) \ldots)x_{n-1})x_n, \
 \ [x_1, \ldots, x_n]_{r(n)} = (x_1(x_2 \ldots (x_{n-1}x_n)\ldots)),
\]
then an $l(n)$-Leibniz derivation ($r(n)$-Leibniz derivation) of $\A$ will be called a~\emph{left (right) Leibniz-derivation of} $\A$ {of order $n$.}
If $d$ is an $f$-derivation for any arrangement $f$ of length $n$, then $d$ will be called a~\emph{Leibniz-derivation of $\A$ of order $n $.} One can see that in our terms a~Leibniz-derivation by Moens is a~right Leibniz-derivation of a~Lie algebra $\A$, and a~Leibniz-derivation by Fialowski-Khudoyberdiev-Omirov is a~left Leibniz-derivation of a~Leibniz algebra $\A $.
It is easy to see that for

\begin{itemize}
    \item[I.] (anti)commutative algebras, the notion of a~left Leibniz-derivation coincides with the notion of a~right Leibniz derivation;

\item[II.] left Leibniz algebras, every right Leibniz-derivation is a~left Leibniz-derivation;

\item[III.] left Zinbiel algebras, every left Leibniz-derivation is a~right Leibniz-derivation;

\item[IV.] associative 
and Lie algebras, the notion of a~Leibniz-derivation coincides with the notion of an $f$-Leibniz derivation for any $f$.
\end{itemize}

 We summarize the main results of the present part in the following theorems.

\begin{theorem}[Theorem 18,~\cite{KP16}]
Let $\A$ be a~Malcev algebra over a~field of characteristic zero.
Then ${\rm Rad}(\A)$ is invariant under all the left Leibniz-derivations of $\A $.
\end{theorem}

\begin{theorem}[Theorem 18,~\cite{KP16}]
Let $\A$ be a~semisimple Malcev algebra over a~field of characteristic 0. Then $\mathfrak{Der}(\A) = {\rm LDer}_l(\A) $.
\end{theorem}

\begin{theorem}[Theorem 28,~\cite{KP16}]
A Malcev algebra over a~field of characteristic zero is nilpotent if and only if it admits an invertible left Leibniz-derivation.
\end{theorem}

\begin{theorem}[Theorem 30,~\cite{KP16}]
A Jordan algebra over a~field of characteristic zero is nilpotent if and only if it admits an invertible left Leibniz--derivation.
\end{theorem}

\begin{theorem}[Theorem 32,~\cite{KP16}]
A $(-1,1)$-algebra over a~field of characteristic zero is nilpotent if and only if it admits an invertible left Leibniz--derivation.
\end{theorem}

\begin{theorem}[Theorem 33,~\cite{KP16}]
A right alternative algebra over a~field of characteristic zero admitting an invertible Leibniz-derivation is right nilpotent.
\end{theorem}

\subsection{Generalized derivations of non-associative algebras}

In 1998, the notion of $\delta$-derivations appeared in the paper by Filippov~\cite{fil1}.
He studied $\delta$-derivations of prime Lie algebras~\cite{fil2}.
After that, $\delta$-derivations of
structurable algebras, $n$-ary algebras, and
Jordan and Lie superalgebras were studied (see~\cite{zus10} and references therein).
The notion of $\frac{1}{2}$-derivation plays an important role in the description of transposed Poisson structures on a~certain Lie algebra~\cite{fkl21}.
The notion of generalized derivations is a~generalization of $\delta$-derivations and $(\alpha,\beta,\gamma)$-derivations.
The most important and systematic research on the generalized derivations algebras of a~Lie algebra and their subalgebras was in a~paper by Leger and Luks~\cite{LL00}.
In their article, they studied the properties of generalized derivation algebras and their subalgebras, for example, the quasiderivation algebras. They determined the structure of algebras of quasiderivations and generalized derivations and proved that the quasiderivation algebra of a~Lie algebra can be embedded into the derivation algebra of a~larger Lie algebra.
Their results were generalized by many authors. For example, Zhang and Zhang~\cite{ZZ10} generalized the above results to the case of Lie superalgebras;
Chen, Ma, Ni and Zhou considered the generalized derivations of color Lie algebras, Hom-Lie superalgebras and Lie triple systems~\cite{cml13}.
Generalized derivations of simple algebras and superalgebras
were investigated in~\cite{ lesha12,lesha14,GP03}.
P\'{e}rez-Izquierdo and Jim\'{e}nez-Gestal used the generalized derivations to study non-associative algebras~\cite{GP08}.
Derivations and generalized derivations of $n$-ary algebras were considered in many papers.
For example, Williams proved that, unlike the case of binary algebras, for any $n \geq 3$ there exists a~non-nilpotent $n$-Lie algebra with invertible derivation~\cite{Williams}.

The main purpose of the results contained in the present part is to generalize the results of Leger and Luks~\cite{LL00} to the case of color $n$-ary algebras and ${\rm Hom}$-algebras.
Particularly, we proved some properties of generalized derivations of color $n$-ary algebras and ${\rm Hom}$-algebras, that the quasiderivation algebra of a~color $n$-ary $\Omega$-algebra (resp., $\Omega$-${\rm Hom}$-algebra) can be embedded into
the derivation algebra of a~large color $n$-ary $\Omega$-algebra (resp., $\Omega$-${\rm Hom}$-algebra)~\cite{kaypopov16,bkp19},
and obtained the classification of all $n$-ary algebras with the property
${\rm End=QDer}$. Let us give the main definitions of this part.

\begin{definition}
A linear mapping $\mathfrak{D}\in {\rm End}(\A)$ is called a~generalized
 derivation of an $n$-ary algebra $\A$ if there exists linear mappings
$\mathfrak{D}',\mathfrak{D}'',\ldots, \mathfrak{D}^{(n-1)},\mathfrak{D}^{(n)} \in {\rm End}(\A)$ such that
\begin{center}
$\sum [x_1, \ldots, \mathfrak{D}^{(i-1)}(x_i), \ldots, x_n]=\mathfrak{D}^{(n)}([x_1, \ldots, x_n]) $.
\end{center}
An $(n+1)$-tuple $(\mathfrak{D}, \mathfrak{D}', \ldots, \mathfrak{D}^{(i-1)}, \ldots, \mathfrak{D}^{(n-1)},\mathfrak{D}^{(n)})$ is called an $(n+1)$-ary derivation.
An $(n+1)$-ary derivation is trivial if it is a~sum of
$(d, \ldots, d)$ and $(\phi_1, \ldots, \phi_{n}, \sum \phi_i) $,
where $d$ is a~derivation and $\phi_i$ are elements from the centroid.
\end{definition}

\begin{definition}
A linear mapping $\mathfrak{D}\in {\rm End}(\A)$ is said to be a~quasiderivation if there exists a~$\mathfrak{D}'\in
{\rm End}(\A)$ such that
\begin{center}
$\sum [x_1, \ldots, \mathfrak{D}(x_i), \ldots, x_n] =\mathfrak{D}^{'}([x_1, \ldots, x_n]) $.
\end{center}
\end{definition}

\begin{definition}
A pair of linear mappings $(d, f)$ satisfying the condition of the previous definition is called a pair of quasiderivations of $\A$.
The image of the projection of ${\rm QDer}(\A)$ onto the first coordinate will be denoted as ${\rm QDer}_{KS}(\A)$ and will be called a~space of quasiderivations in the sense of Kaygorodov and Shestakov~\cite{kay14mz,lesha12},
and the image of the projection of ${\rm QDer}(\A)$ onto the second coordinate will be denoted as ${\rm QDer}_{LL}(\A)$ and will be called a~space of quasiderivations in the sense of Leger and Luks~\cite{LL00}.
\end{definition}

In the following theorems, we summarize the main results of the present part.

\begin{theorem}[Theorem 1,~\cite{kay14st}]
There are no nontrivial ternary derivations of the simple Malcev algebra ${\rm M}_7 $.
\end{theorem}

\begin{theorem}[Theorem 5,~\cite{kay14st}]
There are no nontrivial $4$-ary derivations of the simple ternary Malcev algebra ${\rm M}_8$.
\end{theorem}

\begin{theorem}[Theorem 7,~\cite{kay14mz}]
 Let $\A$ be a~semisimple finite-dimensional $n$-ary Filippov algebra over an algebraically closed field of characteristic zero.
Then, the algebra of generalized derivations of $\A$ is isomorphic to $\oplus \mathfrak{sl}_{n+1} $.
\end{theorem}

To understand the main result of this part, let us give some additional information about $n$-ary algebras (see~\cite{Fil}).

$\bullet$ Up to isomorphism, there is only one $n$-ary anticommutative algebra $A_n$ in dimension $n $. The product of the basis elements $e_1, \dots, e_n$ of $A_n$ is given in the following way:
\[
[e_1,\dots,e_n] = e_1.
\]

$\bullet$ Up to isomorphism, there is only one perfect (i.e., $\A^2 =\A$) $n$-Lie algebra $D_{n+1}$ of dimension $n+1$. The product of the basis elements $e_1, \dots, e_{n+1}$ of $D_{n+1}$ is given in the following way:
\[
[e_1, \dots, \hat{e}_i, \dots, e_{n+1}] = (-1)^{n+i+1}e_i.
\]

$\bullet$ Let $\A$ be an $n$-ary $(n+1)$-dimensional algebra with the basis $e_1, \dots, e_{n+1} $. Let
\begin{equation}\label{ei}
e^i = (-1)^{n+i+1}[e_1,\dots,\hat{e}_i,\dots,e_n], i = 1, \dots, n+1.
\end{equation}
Then, the multiplication in $\A$ is defined by the matrix $B = (\beta_{ij})$ which is given by \[ e^i = \beta_{1i}e_1 + \ldots + \beta_{n+1 i}e_{n+1}, \]
or in terms of matrices as
$
(e^1,\ldots,e^{n+1}) = (e_1,\ldots,e_{n+1})B.
$
The rank of $B$ is equal to the dimension of $\A^2 = [\A,\ldots,\A] $.
An $(n+1)$-dimensional anticommutative algebra with multiplication defined by (\ref{ei}) will be denoted by $\A_B $. It is easy to see that $D_{n+1}$ is an algebra $\A_I$, where $I$ is an identity matrix of order $n+1 $.

Later in the discussion we will also need the description of $1$-dimensional $n$-ary algebras. One can easily see that the multiplication in such algebras is completely defined by an element $\alpha \in \mathbb F$, where $\mathbb F$ is the base field, for it is enough to determine
$[v,v,\ldots,v] = \alpha v$
for any nonzero $v \in \A$ and extend the multiplication linearly. We denote such algebras by $\A_{\alpha}$. It is also easy to see that $\A_{\alpha} \cong \A_{\beta}$ for $\alpha, \beta \neq 0$ if and only if the polynomial $x^{n-1} - \frac{\alpha}{\beta}$ has a~root in $\mathbb F $.

\begin{theorem}[Theorem 5.6,~\cite{kaypopov16}]
Let $\A$ be an $n$-ary algebra such that ${\rm QDer}_{LL}(\A) = End(\A)$. Then

 $\bullet$ If $\A$ is commutative, then either $\A$ has zero product or $\A \cong \A_{\alpha}$ for some $\alpha$ in base field $\mathbb F $. Moreover, if the base field $\mathbb F$ of $\A$ is algebraically closed, then one can define a~binary multiplication $\cdot$ on $\A$ such that $(\A,\cdot) \cong \mathbb F$ and $[x_1,\ldots,x_n] = x_1\cdot\ldots\cdot x_n$, where $[\cdot,\ldots,\cdot]$ is multiplication in $\A$ and $x_1, \dots, x_n \in \A $.

 $\bullet$ If $\A$ is anticommutative, then either $\A$ has zero product or $\A$ is isomorphic to either $A_n$ or $\A_B$ for a~nondegenerate $(n+1) \times (n+1)$-matrix $B$. Moreover, in the last case if $\A$ is a~Filippov algebra and the base field is algebraically closed, then $\A \cong D_{n+1} $.
 \end{theorem}

\begin{theorem}[Theorem 5.9,~\cite{kaypopov16}]
Let $\A$ be an $n$-ary algebra such that ${\rm QDer}_{KS} (\A) = {\rm End}(\A) $. Then

$\bullet$ If $\A$ is commutative, then $\A \cong \A_{\alpha}$ for some $\alpha \in \mathbb F $. Moreover, if the base field $\mathbb F$ of $\A$ is algebraically closed, then one can define a~binary multiplication $\cdot$ in $\A$ such that $(\A,\cdot) \cong \mathbb F$ and $[x_1,\ldots,x_n] = x_1\cdot\ldots\cdot x_n$, where $[\cdot,\ldots,\cdot]$ is multiplication in $\A$ and $x_1, \dots, x_n \in \A $.

$\bullet$ If $\A$ is anticommutative, then $\A \cong \A_B$ for a~nondegenerate $(n+1) \times (n+1)$-matrix $B$ or $\A$ has zero multiplication.
Moreover, if $\A$ is a~Filippov algebra and the base field is algebraically closed, then $\A \cong D_{n+1} $.
\end{theorem}

\section{Poisson type algebras and superalgebras}
The present part is based on the papers
written together with
Bruno Ferreira,
Ivan Shestakov,
Mykola Khrypchenko,
Ualbai Umirbaev
and
Viktor Lopatkin~\cite{fkl21,KK21p,ksu18,K17}.

\subsection{Poisson structures on finitary incidence algebras}

	The systematic study of noncommutative Poisson algebra structures began in the paper by Kubo~\cite{kubo96}.
	He obtained a~description of all the Poisson structures on the full and upper triangular matrix algebras, which was later generalized to prime associative algebras in~\cite{fale98}.
	Namely, it was proved in~\cite{fale98} that any Poisson bracket on a~prime noncommutative associative algebra is the commutator bracket multiplied by an element from the extended centroid of the algebra.
	On the other hand, in his next paper, Kubo studied
	noncommutative Poisson algebra structures on affine Kac-Moody algebras.
	The investigation of Poisson structures on associative algebras continued in some papers by
	Yao, Ye and Zhang~\cite{yyz07}; Mroczyńska, Jaworska-Pastuszak, and Pogorzały~\cite{jpp20,mp18},
	where Poisson structures on finite-dimensional path algebras
	and on canonical algebras were studied.
    In our joint work with Khrypchenko we described Poisson structures on finitary incidence algebras~\cite{KK21p}.

\begin{definition}
	Let $(\A, \cdot)$ be an associative algebra and $\{\cdot, \cdot\}$ an additional multiplication on the same vector space.
	We call $\{\cdot, \cdot\}$ a~Poisson structure on $(\A, \cdot)$ if
	$(\A, \cdot, \{\cdot, \cdot\})$ is a~Poisson algebra.
\end{definition}

    The main tool for the description of Poisson structures on associative algebra is the study of biderivations.
\begin{definition}
	Let $\A$ be an associative algebra and $M$ an $\A$-bimodule.
	A bilinear map $B:\A\times \A \to M$ is called an \textit{antisymmetric biderivation of $\A$ with values in $M$}, if it is
	anticommutative and a~derivation $\A\to M$ with respect to each of its two variables, i.e.
\[
B(xy,z)=B(x,z)y+xB(y,z), \ B(x,x)=0.
\]
For any $\lambda \in C(\A)$, the map
$	B(x,y)=\lambda [x,y],
$	where $[x,y]$ is the commutator $xy-yx$, is an antisymmetric biderivation of $\A$. Such biderivations will be called \textit{inner}.
\end{definition}
For any $\lb\in C(\A)$, the inner biderivation $\{a,b\}=\lb[a,b]$ is clearly a~Poisson structure on $\A$.
Following~\cite{yyz07}, we call such Poisson structures \textit{standard}. There is a~generalization of this notion introduced in the same paper~\cite{yyz07}. Namely, a~Poisson structure on $\A$ is said to be \textit{piecewise standard} if $\A$ decomposes into a~direct sum $\bigoplus_{i=1}^m \A_i$ of indecomposable Lie ideals in such a~way that $\{a,b\}=\lb_i[a,b]$ for all $a\in \A_i$ and $b\in \A$, where $\lb_i\in C(\A)$.

	It is well-known that the \textit{upper triangular matrix algebra} $T_n(R)$ is the incidence algebra $I(C_n,R)$ of a~chain $C_n$ of cardinality $n$. There is a~description of (not necessarily antisymmetric) biderivations of $T_n(R)$ given by Benkovi\v{c} in~\cite[Corollary 4.13]{Benkovic09} (case $n\ge 3$) and~\cite[Proposition 4.16]{Benkovic09} (case $n=2$). A straightforward calculation, based on this description, shows that all the antisymmetric biderivations of $T_n(R)$ are inner. Consequently, we have the following (see also~\cite{kubo96}).

\begin{exm}
		All the Poisson structures on $T_n(R)$ are standard.
	
\end{exm}

	This is not the case for a~general incidence algebra, as the next easy example shows.
	
\begin{exm}\label{2-crown}
		Let $R$ be a~commutative ring and consider $P=\{1,2,3,4\}$ with the following Hasse diagram (a crown).
		
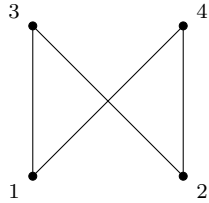
\begin{center}
			
\begin{tikzpicture}[line cap=round,line join=round,>=triangle 45,x=1cm,y=1cm]
			\draw (-1,-1)-- (1,1);
			\draw (1,1)-- (1,-1);
			\draw (1,-1)-- (-1,1);
			\draw (-1,1)-- (-1,-1);
			
\begin{scriptsize}
				\draw [fill=black] (-1,-1) circle (1.5pt);
				\draw[color=black] (-1.25,-1.2) node {$1$};
				\draw [fill=black] (-1,1) circle (1.5pt);
				\draw[color=black] (-1.25,1.2) node {$3$};
				\draw [fill=black] (1,1) circle (1.5pt);
				\draw[color=black] (1.25,1.2) node {$4$};
				\draw [fill=black] (1,-1) circle (1.5pt);
				\draw[color=black] (1.25,-1.2) node {$2$};
			
\end{scriptsize}
			
\end{tikzpicture}
		
\end{center}

\end{exm}

\begin{remark}[Remark 1.3,~\cite{KK21p}]
	 If $R$ is a~field, then the Poisson structure $B$ from \cref{2-crown} is not even piecewise standard.
	
\end{remark}

Let us recall that the definition and some main properties of incidence algebras were considered in Section~\ref{secincidence}.
	For any pair $x\le y$ we shall identify $e_{xy}$ with $1_Re_{xy}\in FI(P,R)$ and denote by $E$ the set $\{e_{xy}\mid x\le y\}\sst FI(P,R)$.
		Denote by $\tilde I(P,R)$ the subalgebra of $FI(P,R)$ generated by $E$. Clearly, $\tl I(P,R)=\textup{span} \{E\}$ as an $R$-module. We will first deal with biderivations defined on $\tilde I(P,R)$. 

\begin{theorem}[Proposition 3.1,~\cite{KK21p}]
		For all $x\le y$ and $u\le v$ we have
		
\begin{center}
		$B(e_{xy},e_{uv})=\lb(e_{xy},e_{uv})[e_{xy},e_{uv}]$
		
\end{center}
		for some map $\lb:E \times E\to R$.
	
\end{theorem}

\begin{theorem}[Proposition 3.5,~\cite{KK21p}]
		Let $\lambda:E\times E\to R$ be a~symmetric map. Then the bilinear map $B$ given by $	B(e_{xy},e_{uv})=\lb(e_{xy},e_{uv})[e_{xy},e_{uv}]$ is an antisymmetric biderivation of $\tilde I(P,R)$ with values in $FI(P,R)$ if and only if
		
\[
	[e_{xy},e_{x'y'}]\ne 0\ \&\ [e_{uv},e_{u'v'}]\ne 0\ \impl\ \lb(e_{xy},e_{x'y'})=\lb(e_{uv},e_{u'v'})
\]
		holds for any chain $C\subseteq P$.
	
\end{theorem}
	
\begin{definition}
		A map $\sigma:(P \times P)_< \to R$ is said to be \textit{constant on chains} if for any chain $C\sst P$ and for all $x<y$, $u<v$ from $C$ one has $\sigma(x,y)=\sigma(u,v)$.
	
\end{definition}

\begin{theorem}[Theorem 3.7,~\cite{KK21p}]
		There is a~one-to-one correspondence between the antisymmetric biderivations $B$ of $\tilde I(P,R)$ with values in $FI(P,R)$ and the maps $\sigma:(P \times P)_< \to R$, which are constant on chains. More precisely,
		
\begin{align}\label{B(f_g)(x_y)=sg(x_y)[f_g](x_y)}
		B(f,g)(x,y)=
\begin{cases}
		0, & x=y,\\
		\sigma(x,y)[f,g](x,y), & x<y,
\end{cases}
\end{align}
		where $f,g\in \tilde I(P,R)$.
	
\end{theorem}

\begin{theorem}[Theorem 4.3,~\cite{KK21p}]
 		The Poisson structures $B$ on $FI(P,R)$ are in a~one-to-one correspondence with the maps $\sigma:(P \times P)_< \to R$ which are constant on chains. The correspondence is given by \cref{B(f_g)(x_y)=sg(x_y)[f_g](x_y)}, in which $f,g\in FI(P,R)$.
	
\end{theorem}

\subsection{Generalized and generic Poisson algebras and superalgebras}
A~generalization of Poisson algebras
known as generic Poisson algebras
was studied in the papers by
Kolesnikov, Makar-Limanov, Shestakov~\cite{KSML}
and Kaygorodov, Shestakov, Umirbaev~\cite{ksu18}.

\begin{definition}
	The triple $(\A, \cdot, \{\cdot, \cdot\})$ is a~generic Poisson algebra if
	$(\A, \cdot)$ is a~commutative associative algebra,
	$(\A,\{\cdot, \cdot\})$ is an anticommutative algebra and
\[
\{a,bc\}=\{a,b\}c+b\{a,c\}.
\]

\end{definition}

In our paper written, together with Shestakov and Umirbaev~\cite{ksu18},
we establish that many results known true for Poisson algebras also hold for generic Poisson algebras. For example,
we study the general properties of generic Poisson modules and of universal multiplicative
enveloping algebras (Section 2,~\cite{ksu18});
we determine the structure of free generic Poisson algebras and of free generic Poisson fields (Section 3,~\cite{ksu18});
we prove an analog of Makar–Limanov–Shestakov’s Theorem, namely that two Poisson
dependent elements in a~free generic Poisson field are polynomial dependent (Section 4,~\cite{ksu18});
we apply the obtained results to the study of automorphisms of the free generic Poisson algebra $GP\{x, y\}$ with two generators and prove that
the groups of automorphisms of the free generic Poisson algebra with two generators,
of the free Poisson algebra with two generators,
of the free associative algebra with two generators
and
of the free commutative associative algebra with two generators are isomorphic (Section 5,~\cite{ksu18}).

Other generalizations of Poisson superalgebras
include generalized Poisson superalgebras and superalgebras of Jordan brackets.
Every Poisson superalgebra is a~superalgebra of Jordan brackets~\cite{Kantor90,kantor92}.
The identities defining unital superalgebras of Jordan brackets
were described in~\cite{KingMcc92}.
Superalgebras of Jordan brackets
are important in the classification of finite-dimensional simple Jordan superalgebras~\cite{MZ01}.
Kac and Cantarini studied linearly compact simple superalgebras of Jordan brackets~\cite{KacKant07},
Kaygorodov and Zhelyabin studied
$\delta$-superderivations of simple superalgebras of Jordan brackets~\cite{ZheKa11},
Zelmanov, Shestakov, and Mart\'{i}nez
studied the relation between Jordan brackets and Poisson brackets~\cite{marsheze01},
and special superalgebras of Jordan brackets were considered in~\cite{ZheKa11}.
Kantor's construction gives interesting relations between
Novikov--Poisson algebras and Jordan superalgebras~\cite{Zakharov}.

A superalgebra
$\A$
is called a~superalgebra of Jordan brackets
whenever its Kantor double is a~Jordan superalgebra~\cite{Kantor90,kantor92}.
As~\cite{KingMcc92} implies,
a superalgebra
$\A$
whose associative-supercommutative multiplication
$\cdot$
and superanticommutative multiplication
$\{ \cdot,\cdot\}$
satisfy
\[
\{a,bc\}=\{a,b\}c+(-1)^{|a||b|}b\{a,c\}-D(a)bc,
\]
{\small \[
\{a,\{b,c\}\}=\{\{a,b\},c\}+(-1)^{|a||b|}\{b,\{a,c\}\}+D(a)\{b,c\}
+(-1)^{|a||bc|}D(b)\{c,a\}+(-1)^{|c||ab|}D(c)\{a,b\},    
\]}

with
$D(a)=\{a,1\}$,
is a~superalgebra of Jordan brackets.
Observe that~%
$D$
is an even derivation of the superalgebra
$(\A,\cdot)$.
If
$D=0$,
then
$(\A,\cdot,\{\cdot,\cdot\})$
is a~Poisson superalgebra~\cite{Kantor90}.
As~\cite{KacKant07} implies,
on the same vector space
we can consider the new multiplication
$\{a,b\}_D= \{a,b\}+\frac{1}{2}(aD(b)-D(a)b)$
and
the new linear mapping
$D_D=\frac{1}{2}D$.
It is a~new superalgebra with the identities
\[
\{a,bc\}_D=\{a,b\}_Dc+(-1)^{|a||b|}b\{a,c\}_D-D_D(a)bc,
\]
\[
\{a,\{b,c\}_D\}_D=\{\{a,b\}_D,c\}_D+(-1)^{|a||b|}\{b,\{a,c\}_D\}_D.
\]
Every superalgebra satisfying the last two identities
is called a~generalized Poisson superalgebra (also known as a~superalgebra of contact brackets).

In our paper, we construct bases for the $n$-generated free unital superalgebra of Jordan brakets
and the $n$-generated free unital generalized Poisson superalgebra in~\cite[Theorems 2 and 5]{K17}.

The celebrated Amitsur--Levitsky theorem states that
the algebra
$M_k(R)$
of
$k \times k$
matrices over a~commutative ring~%
$R$
satisfies the identity
$s_{2k}=0$.
Furthermore,
every associative PI algebra satisfies
$(s_k)^l =0$
by Amitsur's theorem.
Farkas defined customary polynomials
\begin{eqnarray*}\label{tojdest}
g=
\sum_{\sigma \in {\mathbb S}_{m}} c_{\sigma} \{x_{\sigma(1)}, x_{\sigma(2)} \} \ldots \{ x_{\sigma(2i-1)}, x_{\sigma(2i)} \}
\end{eqnarray*}
and proved that
every Poisson PI algebra satisfies some customary identity~\cite{Fr1}.
Farkas' theorem was established for generic Poisson algebras in~\cite{KSML}.

We have found an~analog 
of customary identities
for generalized Poisson algebras and algebras of Jordan brackets~\cite{K17}:
\begin{eqnarray}\label{tojdest}
g_*=
\sum\limits_{i=0}^{[m/2]}\sum_{\sigma \in {\mathbb S}_{m}} c_{\sigma,i} \langle x_{\sigma(1)}, x_{\sigma(2)}\rangle \ldots
\langle x_{\sigma(2i-1)}, x_{\sigma(2i)} \rangle D(x_{\sigma(2i+1)})\ldots D(x_{\sigma(m)}),
\end{eqnarray}
where
\[
\langle x,y\rangle :=\{x,y\} -(D(x)y-xD(y)).
\]

\begin{theorem}[Theorem 13,~\cite{K17}]
If a~unital generalized Poisson algebra~%
$\A$
satisfies a~polynomial identity
$g_0 $,
then~%
$\A$
satisfies a~polynomial identity
$g_*$
of type (\ref{tojdest}).
\end{theorem}

Consequently,
we obtain the following statement.

\begin{theorem}[Corollary 14,~\cite{K17}]
Every PI generalized Poisson algebra
satisfies an identity of the type
{\footnotesize \begin{eqnarray}\label{pi}
 f_{*} =
\sum\limits_{i=0}^{[m/2]}\sum_{\sigma \in \mathbb{S}_{m}} c_{\sigma,i}
 \prod\limits_{k=1}^{i} [x_{\sigma(2k-1)}; x_{\sigma(2k)}; z_{2k-1}; z_{2k}] \cdot
 \prod\limits_{k=1}^{m-2i} \{x_{\sigma(2i+k)}; z_{2i+2k-1}; z_{2i+2k} \} \cdot
 \prod\limits_{k=1}^{2i} z_{2m-2i+k},
\end{eqnarray}}
where
{\small
\[
[u_1; u_2; w_1;w_2]=
\{u_1,u_2\}w_1w_2+\{u_1,w_1w_2\}u_2 + u_1\{w_1w_2, u_2\}+
\sum _{\sigma_1, \sigma_2 \in \mathbb{S}_2} \{u_{\sigma_1(1)}, w_{\sigma_2(1)}\}u_{\sigma_1(2)}w_{\sigma_2(2)},
\]
\[
\{t_1;t_2;t_3 \} = \{t_2t_3,t_1\} - \{t_2,t_1\}t_3-\{t_3,t_1\}t_2.
\]}
\end{theorem}

\begin{theorem}[Theorem 15,~\cite{K17}]
If a~unital algebra of Jordan brackets~%
$\A$
satisfies a~polynomial identity
$g_0$,
then~%
$\A$
satisfies polynomial identities of types (\ref{tojdest})~and~(\ref{pi}).
\end{theorem}

\subsection{Transposed Poisson algebra structures}
Recently, a~dual notion of the Poisson algebra (transposed Poisson algebra) has been introduced in the paper by Bai, Bai, Guo, and Wu~\cite{bai20}. Roughly speaking, this new type of algebra is defined by exchanging the roles of the two binary operations in the Leibniz rule defining the Poisson algebra.
The authors showed that the transposed Poisson algebras defined in this way not only share common properties with Poisson algebras, including the closure under taking tensor products and the Koszul self-duality as an operad but also admit a~rich class of identities. More significantly, a~transposed Poisson algebra naturally arises from a~Novikov-Poisson algebra by taking the commutator Lie algebra of the Novikov algebra. Consequently, we find a~parallelism between the classical construction of a~Poisson algebra from a~commutative associative algebra with a~pair of commuting derivations and the construction of a~transposed Poisson algebra when there is some non-zero derivation. Furthermore, the transposed Poisson algebras also capture the algebraic structure when the commutator is taken in pre-Lie Poisson algebras and two other algebras of Poisson type.

\begin{definition}\label{tpa}
Let $\A$ be a~vector space equipped with two nonzero bilinear operations $\cdot$ and $[\cdot,\cdot] $.
The triple $(\A,\cdot,[\cdot,\cdot])$ is called a~transposed Poisson algebra if $(\A,\cdot)$ is a~commutative associative algebra and
$(\A,[\cdot,\cdot])$ is a~Lie algebra that satisfies the following compatibility condition
\begin{equation*}
2z\cdot [x,y]=[z\cdot x,y]+[x,z\cdot y].
\end{equation*}
\end{definition}

The study of $\delta$-derivations of Lie algebras was initiated by Filippov in 1998~\cite{fil1,fil2}.
Let us define $\frac{1}{2}$-derivations which will be very useful for our work.

\begin{definition}\label{12der}
Let $(\A, [\cdot,\cdot])$ be an algebra with multiplication $[\cdot,\cdot]$ and $\varphi$ be a~linear map.
Then $\varphi$ is a~$\frac{1}{2}$-derivation if it satisfies
$\varphi[x,y]= \frac{1}{2} \left([\varphi(x),y]+ [x, \varphi(y)]\right) $.
\end{definition}

Summarizing definitions~\ref{tpa} and~\ref{12der}, we have the following key statement.

\begin{theorem}[Lemma 7,~\cite{fkl21}]
Let $(\A,\cdot,[\cdot,\cdot])$ be a~transposed Poisson algebra
and $z$ an arbitrary element from $\A $.
Then the right multiplication $R_z$ in the commutative associative algebra $(\A,\cdot)$ gives a~$\frac{1}{2}$-derivation of the Lie algebra $(\A, [\cdot,\cdot]) $.
\end{theorem}

In our joint work with Ferreira and Lopatkin~\cite{fkl21}, we found a~way to describe all transposed Poisson algebra structures with a~certain Lie algebra.
Our main tools are the description of the space of $\frac{1}{2}$-derivations of Lie algebras
and the connection between the space of $\frac{1}{2}$-derivations of the Lie part of a~transposed Poisson algebra and the space of right multiplications of the associative part of this transposed Poisson algebra. Namely, every right \textit{associative} multiplication is a~\textit{Lie} $\frac{1}{2}$-derivation.
Using the known description of $\delta$-derivations of semisimple finite-dimensional Lie algebras,
we have found that there are no transposed Poisson algebras with a~semisimple finite-dimensional Lie part.
In the case of simple infinite-dimensional algebras, we have a~different situation: it has been proved that the Witt algebra admits many nontrivial structures of transposed Poisson algebras.
Later we studied structures of transposed Poisson algebras defined on one of the most interesting generalizations of the Witt algebra.
Namely, we considered $\frac{1}{2}$-derivations of the algebra ${\mathcal W}(a,b)$ and proved that the algebra
${\mathcal W}(a,b)$ does not admit structures of transposed Poisson algebras if and only if $b\neq -1 $.
All the transposed Poisson algebra structures defined on ${\mathcal W}(a,-1)$ have been described.
In the following section, we considered another example of an algebra related to the Witt algebra.
We proved that there are no transposed Poisson algebras defined on the Virasoro algebra.
As for some corollaries, we proved that there are no transposed Poisson algebra structures defined on the $N=1$ and $N=2$ superconformal algebras. The rest of the paper is dedicated to a~classification of $\frac{1}{2}$-derivations and constructions of transposed Poisson algebras defined on the
thin Lie algebra and the solvable Lie algebra with the abelian nilpotent radical of codimension $1 $.
Other interesting examples of transposed Poisson algebras constructed in
a series of papers together with Khrypchenko~\cite{kkk1,kkk2,kkk3,kkk4}.

\section{$n$-Ary algebras}
The present part is based on the papers
written together with
Alexandre Pozhidaev,
Antonio Jesús Calderón,
Elisabete Barreiro,
José María Sánchez,
Paulo Saraiva,
and
Yury Popov~\cite{bks19,bcks19,cks19,KP19,KPS19}.

\subsection{Ternary Jordan algebras}
 Based on the relation between the notions of Lie triple systems and Jordan algebras, we introduce the $n$-ary Jordan algebras~\cite{KPS19},
an $n$-ary generalization of Jordan algebras obtained via the generalization of the following property $\left[ R_{x},R_{y}\right] \in {\mathfrak Der}\left(\A\right) $,
where $\A$ is an $n$-ary algebra. As we can see, the present generalization of algebras is an $n$-ary generalization of almost Jordan and commutative $\mathfrak{CD}$-algebras.

In~\cite{KPS19}, we give the first example of ternary Jordan algebras defined in this way and study some properties of this algebra.
Namely, let $\mathbb{V}$ be an $n$-dimensional vector space over a~field $\mathbb{F}$, and denote the bilinear form by $\left(\cdot,\cdot\right)$. Consider the following ternary multiplication defined on $\mathbb{V}$:
\begin{equation*}
\llbracket x,y,z \rrbracket =\left(y,z\right) x+\left(x,z\right) y+\left(
x,y\right) z.
\end{equation*}
Denote the obtained ternary algebra by $\mathbb{A}$.
The following three theorems give us a~characterization of $\mathbb{A}$.

\begin{theorem}[Theorem 5,~\cite{KPS19}]
$\mathbb{A}$ is a~ternary Jordan algebra. \label{THM 1}
\end{theorem}

\begin{theorem}[Theorem 6,~\cite{KPS19}]
The ternary Jordan algebra $\mathbb{A}$ is simple, except if $\dim\ \mathbb{V}=2$ and $char \left(\mathbb{F}\right) = 2$.
\end{theorem}

\begin{theorem}[Theorem 11,~\cite{KPS19}]\label{der}
$\mathfrak{Der}\left(\mathbb{A}\right) ={\rm Inder}\left(\mathbb{A}\right)=\mathfrak{so}(n) $.
\end{theorem}

\begin{remark}[Remark 12,~\cite{KPS19}]
In $1955$ Jacobson proved that if a~finite-dimensional Lie algebra over a~field of characteristic zero has an invertible derivation,
then it is a~nilpotent algebra~\cite{Jac}.
The same result was proved for Jordan algebras~\cite{KP16}, but as we can see from Theorem~\ref{der},
the Theorem of Jacobson is not true for ternary Jordan algebras.
We can take the ternary Jordan algebra $\mathbb{A}$ (as in Theorem~\ref{der}) with dimension $4$ and
consider the map defined by the matrix $\sum_{1\leq i<j\le 4} (e_{ij}-e_{ji}) $.
As follows, there is a~simple ternary Jordan algebra with an invertible derivation.
\end{remark}

Later in~\cite{KPS19}, we study a~ternary example of these algebras. Finally, based on the construction of a~family of ternary algebras defined using the Cayley --- Dickson algebras, we present an example of a~ternary $\mathfrak{CD}$-algebra ($n$-ary $\mathfrak{CD}$-algebras are the $n$-ary version of (not necessarily commutative) $\mathfrak{CD}$-algebras) in~\cite[Section 5]{KPS19}.

Let $\mathcal{U}_{2}$ be the generalized quaternions
and
$\mathcal{U}_{3}$ be the generalized octonions.
Define on $\mathcal{U}_{t}$, for $t=2,3$, the ternary
multiplication:
\[
\llbracket x,y,z\rrbracket =\left(x\overline{y}\right) z,
\]
and take
$\mathcal{D}_{t}=\left(\mathcal{U}_{t},\llbracket \cdot, \cdot, \cdot\rrbracket\right) $.

This ternary multiplication is not totally commutative,
so these algebras are not ternary Jordan algebras.

\begin{theorem}[Theorem 15,~\cite{KPS19}]
$\mathcal{D}_{2}$ is a~simple ternary $\mathfrak{CD}$-algebra.
\end{theorem}

\begin{theorem}[Lemma 15,~\cite{KPS19}]
  $\mathcal{D}_{3}$ is not a~ternary $\mathfrak{CD}$-algebra.
\end{theorem}

\subsection{Split regular ${\rm Hom}$-Leibniz color $3$-algebras}

The study of ${\rm Hom}$-structures began in the paper by Hartwig, Larsson, and Silvestrov~\cite{hom}.
The notion of ${\rm Hom}$-Lie triple systems was introduced in~\cite{[20]}.
In the paper~\cite{[28]} Yau gave a~general method for constructing ${\rm Hom}$-type algebras starting from usual algebras and a~twisting self-map.

 In our joint work with Popov~\cite{KP19},
 we study the structure of split regular ${\rm Hom}$-Leibniz $3$-algebras of arbitrary dimension and over an arbitrary base field ${\mathbb F}$.
  Split structures first appeared in the classical theory of (finite-dimensional) Lie algebras, but have been extended to more general settings
  like, for example,
  Leibniz algebras~\cite{AJC},
  Poisson algebras,
  Leibniz superalgebras,
  regular ${\rm Hom}$-Lie algebras,
  regular ${\rm Hom}$-Lie superalgebras,
  regular ${\rm Hom}$-Lie color algebras,
  regular ${\rm Hom}$-Poisson algebras~\cite{AraCM2},
  regular ${\rm Hom}$-Leibniz algebras,
  regular ${\rm BiHom}$-Lie algebras~\cite{AJCS2}, and
  regular ${\rm BiHom}$-Lie superalgebras,
  among many others.
    As for the study of split ternary structures, see~\cite{AJF2} for Lie triple systems, twisted inner derivation triple systems, Lie $3$-algebras~\cite{AJF2}, Leibniz $3$-algebras~\cite{3L}, and for Leibniz triple systems.

\begin{definition}
{
Let ${\mathbb F}$ be a~field and $\mathbb{G}$ be an abelian group.
A map $\epsilon: \mathbb{G} \times \mathbb{G} \rightarrow {\mathbb F}^{\times}$ is called a~bicharacter on $\mathbb{G}$ if the
following relations hold for all $f,g,h \in \mathbb{G}:$

(1) $\epsilon(f,g+h)=\epsilon(f,g)\epsilon (f,h) $;

(2) $\epsilon(g+h,f)=\epsilon(g,f)\epsilon(h,f) $;

(3) $\epsilon(g,h)\epsilon(h,g)=1 $.}
\end{definition}

\begin{definition}
A $\mathbb{G}$-graded color $n$-ary algebra $T$ is a~vector space $T=\bigoplus_{g\in \mathbb{G}} T_g$ with an $n$-linear map
$[\cdot, \ldots, \cdot ]: T \times \ldots \times T \rightarrow T$ satisfying
\[
[T_{\theta_1}, \ldots, T_{\theta_n}] \subseteq T_{\theta_1+ \ldots +\theta_n}, \ \theta_i \in \mathbb{G}.
\]
\end{definition}

\begin{definition}\label{colorleibter}
A regular ${\rm Hom}$-Leibniz color $3$-algebra $(T,[\cdot,\cdot,\cdot], \epsilon, \phi)$ is a~$\mathbb{G}$-graded vector space $T$ with a~bicharacter $\epsilon$, an even trilinear map $[\cdot,\cdot,\cdot]$ and an even automorphism $\phi$ satisfying
\[
[\phi(x_1),\phi(x_2),[y_1,y_2,y_3]]=
[[x_1,x_2,y_1],\phi(y_2),\phi(y_3)] +
\]
\[
\epsilon(x_1 + x_2,y_1)[\phi(y_1),[x_1,x_2,y_2],\phi(y_3)] + \epsilon(x_1+x_2,y_1+y_2)[\phi(y_1),\phi(y_2),[x_1,x_2,y_3]].
\]

Additionally, if the identities
\[
[x_1,x_2,x_3]=-\epsilon(x_1,x_2)[x_2,x_1,x_3],\ [x_1,x_2,x_3]=-\epsilon(x_2,x_3)[x_1,x_3,x_2]
\]
hold in $T$, then $T$ is called a~regular ${\rm Hom}$-Lie color $3$-algebra, and if the identities
\begin{gather*}[x_1,x_2,x_3]
=-\epsilon(x_1,x_2)[x_2,x_1,x_3],\\ \epsilon(x_3,x_1)[x_1,x_2,x_3]+ \epsilon(x_1,x_2)[x_2,x_3,x_1]+ \epsilon(x_2,x_3)[x_3,x_1,x_2]=0
\end{gather*}
hold in $T $, then $T$ is called a~regular ${\rm Hom}$-Lie triple color system.

\begin{definition}\label{defembed}
Let $T = (T,[\cdot,\cdot,\cdot],\epsilon,\phi)$ be a ${\rm Hom}$-Leibniz color $3$-algebra.
Consider the space $\mathfrak{L} = {\rm span}\{\ad(x,y): x, y \in T\} $, where $\ad (x,y)(z):=[x,y,z] $.
The \textit{ standard embedding} of a~regular ${\rm Hom}$-Leibniz color $3$-algebra $(T,[\cdot,\cdot,\cdot],\epsilon,\phi)$ is a~color
$2$-graded algebra with an automorphism $(\A,\epsilon,\Phi)$, where ${\A}^{0}:= {\mathfrak L}$, ${\A}^{1}:= T$, the product is given by
\begin{equation*}
(\ad(x,y),z) \cdot (\ad(u,v),w) =
\end{equation*}
\begin{equation*}
(\ad([x,y,u]^{\phi^{-1}},v) + \epsilon(x+y,u)\ad(u, [x,y,v]^{\phi^{-1}})+ \ad(z,w) \, \ [x,y,w]^{\phi^{-1}}-\epsilon(z, u+v)[u,v,z]^{\phi^{-1}}),
\end{equation*}

the automorphism $\Phi$ by
\begin{equation*}
\Phi:
\begin{cases}
x &\mapsto \phi(x),\\
\ad(x,y) &\mapsto \phi\ad(x,y)\phi^{-1} = \ad(\phi(x),\phi(y)),
\end{cases}
\text{ for } x, y \in T,
\end{equation*}
and the $\mathbb{G}$-grading is induced by the $\mathbb{G}$-gradings of $T$ and $\mathfrak{L} $.
\end{definition}

\begin{definition}\label{root_spaces}
Let $T = (T,[\cdot,\cdot,\cdot],\epsilon,\phi)$ be a ${\rm Hom}$-Leibniz color $3$-algebra (or a~Leibniz color $3$-algebra with an automorphism) and let ${\A} = ({\mathfrak{L}}\oplus T,\cdot,\epsilon,\Phi)$ be its standard embedding. Let $H$ be a~maximal abelian subalgebra (shortly MASA) of $\mathfrak{L}_0 $, the zeroth $\mathbb G$-component of $\mathfrak{L} $. The root space of $T$ with respect to $H$ associated to a~linear functional $\alpha \in H^{*}$ is the subspace
\[
T_{\alpha }:= \{v\in { T}:h\cdot v = \alpha(h) v \, \mbox{ for any } \, h\in H\}.
\]
The elements $\alpha \in H^{*}$ such that $T_{\alpha }\neq 0$ are called roots of $T$ (with respect to $H$), and we
 {write} $\Lambda^T :=\{\alpha \in H^{*}\backslash \{0\}: T_{\alpha }\neq 0\}$.
Analogously, by $\Lambda^{{\mathfrak{L}}}$ we denote the set of all non-zero $\alpha \in H^{*}$ such that ${\mathfrak{L}}_{\alpha } \neq 0$, where
\[
{\mathfrak{L}}_{\alpha }:=\{e\in {\mathfrak{L}}:[h,e]=\alpha (h)e \, \mbox{ for any } \, h\in H\}
\]
are the root subspaces of $\mathfrak{L}$ with respect to $T $.
\end{definition}

\end{definition}

\begin{definition}
Let $T$ be a~${\rm Hom}$-Leibniz color $3$-algebra or a~Leibniz color $3$-algebra with an automorphism and $\A = \mathfrak{L}\oplus T$ its standard embedding. Then $T$ is said to be a~split algebra if there exists a~MASA $H$ of ${\mathfrak{L}_0}$ such that
\begin{equation}\label{rootdeco}
T = T_{0}\oplus \left (\bigoplus _{\alpha \in \Lambda^T}T_{\alpha
}\right).
\end{equation}
in the sense of Definition~\ref{root_spaces}. The set $\Lambda^T$ is called the root system of $T$. We refer to the decomposition~\eqref{rootdeco} as the root spaces decomposition of $T$.
\end{definition}

\begin{definition}
The root system $\Lambda^T$ ($\Lambda^{\mathfrak L}$) is called symmetric if $\Lambda^T = -\Lambda^T$ ($\Lambda^{\mathfrak L} = -\Lambda^{\mathfrak L}$), where for $\emptyset \neq \Upsilon \subset H^*$ the set $-\Upsilon$ is just $\{-\alpha: \alpha \in \Upsilon\}$.
\end{definition}

\begin{theorem}[Theorem 3.16,~\cite{KP19}]
Let $(T,\phi)$ be a~split ${\rm Hom}$-Leibniz color $3$-algebra with multiplication algebra $\mathfrak{L} $. Suppose that the root systems $\Lambda^T$ and $\Lambda^{\mathfrak{L}}$ are symmetric. Then there exists an equivalence relation $\sim$ on $\Lambda^T$ and a~subspace $\mathcal{U} \subseteq T_0$ such that $T = \mathcal{U} \oplus \sum_{[\alpha] \in \Lambda^T /\sim} T_{[\alpha]}$.
Moreover,
\[
[T, T_{[\alpha]},T_{[\beta]}] + [T_{[\alpha]},T,T_{[\beta]}] + [T_{[\alpha]},T_{[\beta]}, T]= 0
\]
whenever $[\alpha] \neq [\beta]$.
\end{theorem}

\subsection{$n$-Ary algebras with a~multiplicative type basis}

In the literature, it is usual to describe an algebra by exhibiting a~multiplicative table among the elements of a~fixed basis. There exist many classical examples of algebras admitting multiplicative bases in the setting of several varieties as associative algebras, Lie algebras, Malcev algebras, Leibniz algebras, ${\rm Hom}$-Lie algebras, etc. For instance, in the class of associative algebras, we have the classes of full matrix algebras, group algebras, quiver algebras, etc. In the class of Lie algebras, we can consider the semisimple finite-dimensional Lie algebras over algebraically closed fields of characteristic $0$,
the Heisenberg algebras, the twisted Heisenberg algebras, and so on.
In the class of Leibniz algebras, we have the classes of (complex) finite-dimensional naturally graded filiform Leibniz algebras and $n$-dimensional filiform graded filiform Leibniz algebras of length $n-1$. By looking at the multiplication table of the non-Lie Malcev algebra ${\rm M}_7$ (a $7$-dimensional simple algebra), we have another example of an algebra with a~multiplicative basis. For Zinbiel algebras we have, for instance, that any complex $n$-dimensional null-filiform Zinbiel algebra admits a~multiplicative basis.
We can also mention the simple finite-dimensional Filippov algebras in~\cite{Fil} as examples of $n$-ary algebras admitting a~multiplicative basis ($n=3,4,\dots$).
Many infinite-dimensional simple Filippov algebras admitting a~multiplicative basis were constructed in papers by Bai, Ding, Pozhidaev, etc.
The concept of multiplicative bases appears in a~natural way in the study of different physical problems. One may expect that there are problems that are naturally and more simply formulated by exploiting multiplicative bases.
 There are many works about the study of algebras with a~multiplicative basis\cite{m6,m5,m3,m1,Yo_n_algebras,Yo_modules,mm,cuasi,vitya03,vitya00}.
 Calderón introduced the notion of quasi-multiplicative bases as a~generalization of the notion of multiplicative bases in~\cite{cuasi} and described the structure of associative algebras with a~quasi-multiplicative basis.

\begin{definition}\label{11}\rm
A basis of homogeneous elements $\mathbb{B} = \{e_i\}_{i \in I}$ of a~graded $n$-ary algebra $\A$ is multiplicative if for any $i_1,\dots,i_n \in I$ we have $[ e_{i_1},\dots,e_{i_n} ] \in \mathbb{F}e_j$ for some $j \in I$.
\end{definition}

In the paper~\cite{Yo_modules} of
Calderón, Navarro Izquierdo, and Sánchez Delgado, the authors considered a~generalization of the multiplicative basis of an algebra.
Let $\mathbb V$ and $\mathbb W$ be two vector spaces over a~base field $\mathbb F$.
It is said that $\mathbb V$ is a~module over $\mathbb W$ if it is endowed with a~bilinear map $\mathbb V \times \mathbb W \to \mathbb V$. A basis $\mathbb B =\{v_i\}_{i \in I}$ of $\mathbb V$ is called multiplicative with respect to another basis $\mathbb B'=\{w_j\}_{j\in J}$ of $\mathbb W$ if for any $i\in I,j \in J$ we have $v_iw_j \in \mathbb{F} v_k$ for some $k\in I$.
They showed that if $\mathbb V$ admits a~multiplicative basis in the above sense, then it decomposes as the direct sum $\mathbb V=\oplus V_{\alpha}$ of well-described submodules, admitting each one a~multiplicative basis. Also, under a~mild condition, the minimality of $\mathbb V$ is characterized in terms of the multiplicative basis and it is shown that the above direct sum is using the family of its minimal submodules, admitting each one a~multiplicative basis.

In our joint paper with Barreiro and Sánchez~\cite{bks19}, we generalize the cited results to $n$-ary case. We study the structure of certain $k$-modules $\mathbb{V}$ over linear spaces $\mathbb{W}$ without restrictions neither on the dimensions of $\mathbb{V}$ and $\mathbb{W}$ nor on the base field $\mathbb F$. A basis $\mathbb{B} = \{v_i\}_{i\in I}$ of $\mathbb{V}$ is called multiplicative with respect to another basis $\mathbb{B}' = \{w_j\}_{j \in J}$ of $\mathbb{W}$ if for any $\sigma \in \mathbb S_n $, $i_1,\dots,i_k \in I$ and $j_{k+1},\dots, j_n \in J$ we have \[ [v_{i_1},\dots, v_{i_k}, w_{j_{k+1}}, \dots, w_{j_n}]_{\sigma} \in \mathbb{F}v_{r_{\sigma}} \] for some $r_{\sigma} \in I$. We show that if $\mathbb{V}$ admits a~multiplicative basis then it decomposes as the direct sum $\mathbb{V} = \bigoplus_{\alpha} V_{\alpha}$ of well-described $k$-submodules $V_{\alpha}$ each one admitting a~multiplicative basis. Also, the minimality of $\mathbb{V}$ is characterized in terms of the multiplicative basis and it is shown that the above direct sum is using the family of its minimal $k$-submodules, admitting each one a~multiplicative basis. Finally, we study an application of $k$-modules with a~multiplicative basis over an arbitrary $n$-ary algebra with a~multiplicative basis.

We continue our study of $n$-ary linear maps on a~vector space in a~joint paper together with Calderón and Saraiva~\cite{cks19}.
Let $\mathbb V$ be an arbitrary linear space and $f:\mathbb V \times \ldots \times \mathbb V \to \mathbb V$ an $n$-linear map. It is proved that, for each choice of a~basis $\mathbb{B}$ of $\mathbb V$, the $n$-linear map $f$ induces a~(nontrivial) decomposition $\mathbb V= \oplus V_j$ as a~direct sum of linear subspaces of $\mathbb V$, with respect to $\mathbb{B}$. It is shown that this decomposition is $f$-orthogonal in the sense that $f(\mathbb V, \ldots, V_j, \ldots, V_k, \ldots, \mathbb V) =0$ when $j \neq k$, and in such a~way that any $V_j$ is strongly $f$-invariant, meaning that $f(\mathbb V, \ldots, V_j, \ldots, \mathbb V) \subset V_j $.
We deduce a~sufficient condition for the existence of an isomorphism between two different decompositions of $\mathbb V$ induced by an $n$-linear map $f$, concerning two different bases of $\mathbb V$. We also characterize the $f$-simplicity -- an analog of the usual simplicity in the framework of $n$-linear maps -- of any linear subspace $V_j$ of a~certain decomposition induced by $f$. Finally, an application to the structure theory of arbitrary $n$-ary algebras is provided. This work is a close generalization of the results obtained by Calderón (2018)~\cite{Yo4}.

The rest of our results in the present topic are dedicated to a~generalization of the paper~\cite{cuasi} of Calderón.
Namely, in a~joint work together with
 Barreiro, Calderón, and Sánchez~\cite{bcks19}, we obtained an analog of Calderón´s results in the class of color generalized Lie type $n$-ary algebras.

\begin{definition}
A graded $n$-ary algebra $\A$ admits a~\textit{quasi-multiplicative basis} if $\A = \mathbb V \oplus \mathbb W$, with $\mathbb V$ and $0 \neq \mathbb W$ graded linear subspaces in such a~way that there exists a~basis of homogeneous elements $\mathbb{B} = \{e_i\}_{i \in I}$ of $\mathbb W$ satisfying:
\begin{itemize}
\item[1.] For $i_1,\dots,i_n \in I$ we have either $[ e_{i_1}, \dots, e_{i_n} ] \in \mathbb{F}e_j$ for some $j \in I$ or $[ e_{i_1},\dots, e_{i_n} ] \in \mathbb V$.
\item[2.] Given $0 < k < n $, for $i_1,\dots,i_k \in I$ and $\sigma \in \mathbb S_n$ we have $[ e_{i_1}, \dots, e_{i_k}, \mathbb V, \dots, \mathbb V ]_{\sigma} \subset \mathbb{F}e_{j_{\sigma}}$ for some $j_{\sigma} \in I$.
\item[3.] We have either $[ \mathbb V, \dots, \mathbb V ] \subset \mathbb{F}e_j$ for some $j \in I$ or $[ \mathbb V, \dots,\mathbb V ] \subset\mathbb V$.
\end{itemize}
\end{definition}

\begin{definition}
\rm
An $n$-ary algebra $(\A, [ \cdot, \dots, \cdot ])$ is called a~\textit{generalized Lie-type algebra} if it satisfies the following $n$ identities:
\begin{eqnarray*}
\begin{split}
& [ y_1, \ldots, \underbrace{[ x_1, \dots, x_n ]}_{\mbox{pos } k}, \dots, y_{n-1} ] =\\
& \quad \quad \quad \sum\limits_{{\tiny
\begin{array}
{l}
  1 \leq i,j \leq n \\
  \sigma_1 \in \mathbb S_n\\
  \sigma_2 \in \mathbb S_{n-1} \\
\end{array}
}}
\alpha_{i,j,k}^{\sigma_1,\sigma_2} [ x_{\sigma_1(1)}, \dots, x_{\sigma_1(i-1)}, [ y_{\sigma_2(1)}, \dots, \underbrace{x_{\sigma_1(i)}}_{\mbox{pos } j}, \dots, y_{\sigma_2(n-1)} ], x_{\sigma_1(i+1)}, \dots, x_{\sigma_1(n)} ],
\end{split}
\end{eqnarray*}
for $k=1, \ldots, n $, being $\alpha_{i,j,k}^{\sigma_1,\sigma_2} \in \mathbb{F}$, and where \mbox{pos} $j$ means that the element $x_{\sigma_1(i)}$ is in the position $j$ in the inside $n$-product.
\end{definition}

 Observe that we can obtain several binary and $n$-ary algebras depending of the values of $\alpha_{i,j,k}^{\sigma_1,\sigma_2}$:
\begin{itemize}
\item[$\bullet$] Lie algebras, Leibniz algebras, Novikov algebras,
associative algebras, alternative algebras, bicommutative algebras, commutative pre-Lie algebras, etc.;
\item[$\bullet$] $n$-Lie (Filippov) algebras, commutative Leibniz $n$-ary algebras,
totally associative-commutative $n$-ary algebras, etc.
\end{itemize}

\begin{theorem}[Theorem 21,~\cite{bcks19}]
A color generalized Lie type algebra $\A = \mathbb V \oplus \mathbb W$ admitting a~quasi-multiplicative basis of $\mathbb W \neq 0$ decomposes as
\[
\A = {\mathcal U} \oplus \Bigl(\sum\limits_{[i] \in I/\sim}{\frak J}_{[i]} \Bigr),
\]
where ${\mathcal U}$ is a~linear complement of $\sum_{[i] \in I/\sim}\mathbb V_{[i]}$ in $\mathbb V$ and any ${\frak J}_{[i]}$ is one of the color generalized Lie-type ideals, admitting a~quasi-multiplicative basis. Furthermore
\[
[ {\frak J}_{[i]}, {\frak J}_{[h]}, \A, \ldots, \A ]_{\sigma} = 0
\]
whenever $[i] \neq [h]$.
\end{theorem}


\EditInfo{June 2, 2023}{August 1, 2023}{Ana Cristina Moreira Freitas, Carlos Florentino, and   Diogo Oliveira e Silva }

\end{document}